\documentclass[11pt]{article}
\usepackage[a4paper,margin=2.7cm]{geometry}
\usepackage[T1]{fontenc}
\usepackage[utf8]{inputenc}
\usepackage{lmodern}
\usepackage{amsmath,amssymb,amsthm,mathtools}
\usepackage{enumitem,booktabs,longtable,array,microtype}
\usepackage{xcolor}
\newif\ifshowchanges
\showchangestrue

\definecolor{revisiongreen}{RGB}{0,0,0}
\usepackage{tikz}
\usetikzlibrary{arrows.meta}
\usepackage[hidelinks]{hyperref}
\numberwithin{equation}{section}
\newtheorem{theorem}{Theorem}[section]
\newtheorem{proposition}[theorem]{Proposition}
\newtheorem{lemma}[theorem]{Lemma}
\newtheorem{corollary}[theorem]{Corollary}
\theoremstyle{definition}
\newtheorem{definition}[theorem]{Definition}
\newtheorem{assumption}[theorem]{Assumption}
\theoremstyle{remark}
\newtheorem{remark}[theorem]{Remark}
\DeclareMathOperator{\dev}{dev}
\DeclareMathOperator{\tr}{tr}

\DeclareMathOperator{\Div}{div}
\DeclareMathOperator{\supp}{supp}
\newcommand{\R}{\mathbb R}
\newcommand{\Ssym}{\mathcal S^3}
\newcommand{\Sdev}{\mathcal S^3_{\rm dev}}
\newcommand{\M}{\mathcal M}
\newcommand{\eps}{\varepsilon}
\newcommand{\Id}{\mathbb I}
\newcommand{\C}{\mathbb C}
\newcommand{\TT}{\mathbb T}
\newcommand{\weakstar}{\overset{*}{\rightharpoonup}}
\newcommand{\dd}{\,\mathrm d}
\newcommand{\Hdiv}{H(\operatorname{div};\Omega)}
\newcommand{\Hzdiv}{H_0(\operatorname{div};\Omega)}

\newcommand{\norm}[1]{\lVert#1\rVert}
\newcommand{\dual}[2]{\langle#1,#2\rangle}
\newcommand{\qbar}{\overline Q}
\newcommand{\ip}[2]{\left\langle #1,#2\right\rangle}
\setlist[enumerate]{label=\textup{(\roman*)},leftmargin=*,itemsep=4pt}
\title{Dynamic Prandtl--Reuss Thermoplasticity with Degenerating Yield Constraints}
\author{Sebastian Owczarek\\
{\small Faculty of Mathematics and Information Science, Warsaw University of Technology,}\\
{\small Koszykowa 75, 00-662 Warsaw, Poland}\\
{\small\texttt{sebastian.owczarek@pw.edu.pl}}\\
{\small ORCID: \href{https://orcid.org/0000-0001-8860-0243}{0000-0001-8860-0243}}}
\date{}
\begin{document}
\maketitle

\begin{abstract}
We study a three-dimensional dynamic thermoplastic system with a temperature-dependent convex yield set whose deviatoric section may collapse completely at a finite temperature. The underlying constitutive model is rate-independent Prandtl--Reuss perfect plasticity, without Kelvin--Voigt viscosity, hardening, or a safe-load condition. We construct global approximation-generated relaxed entropy--energy states by a one-parameter Temam--Norton--Hoff--Galerkin scheme with an undiscretized heat equation and untruncated plastic heat production. The logarithmic temperature is strongly compact in $L^2(0,T;L^q(\Omega))$ for every $q<6$ and in $L^p((0,T)\times\Omega)$ for every $p<8/3$. An entropy-weighted work estimate yields convergence of the overstress at rate $O(\sqrt{\log n/n})$ in every fixed finite Lebesgue space and identifies stress admissibility. The constitutive limit comprises a distributional plastic balance, finite support-weighted measures, and compact-fibre lifts carrying exact normality and maximum plastic work. A locally finite lift on the positive-radius fibres unifies the level-restricted limits at every continuity level. At the macroscopic level, however, the available compactness does not identify a canonical stress--plastic-flow pairing. Moreover, an explicit smooth state satisfying the intrinsic relaxed balance conditions shows that these conditions do not imply macroscopic Prandtl--Reuss normality. In the fully collapsed regime, an energy--entropy rigidity argument determines the evolution and yields approximation-generated examples in which the plastic flow need not be a finite measure or the displacement need not retain its full boundary trace. Thermal concentration and plastic entropy production remain explicit throughout.
\end{abstract}



\section{Introduction}
\label{sec:introduction}
Temperature affects both the mechanical response of a solid and its resistance to irreversible deformation. In a material undergoing thermal softening, plastic dissipation produces heat, while the resulting increase in temperature may reduce the yield strength and thereby alter the subsequent plastic evolution. An extreme, but mathematically natural, idealization is the complete loss of deviatoric strength above a critical temperature.

The mathematical analysis of such thermomechanical systems is substantially more delicate than that of the corresponding isothermal models. In the absence of Kelvin--Voigt viscosity and hardening, the mechanical energy controls the stress and the velocity, whereas the heat equation is driven by a source of essentially $L^1$-type. At the same time, degeneration of the yield set destroys the uniform coercivity of the plastic dissipation and leaves no global $L^1$-control of the plastic rate.

We study this regime in a three-dimensional small-strain dynamic model. Let $\Omega\subset\mathbb R^3$ be a bounded domain and let $T>0$. On the spacetime cylinder $Q=(0,T)\times\Omega,$ the unknowns are the displacement $u:Q\to\mathbb R^3$, the elastic stress $\TT:Q\to\mathcal S^3$, the absolute temperature $\theta:Q\to(0,\infty)$, and the deviatoric plastic strain rate $\xi:Q\to\mathcal S^3_{\rm dev}$. Formally, they satisfy
\begin{equation}\label{eq:intro-model}
\begin{aligned}
u_{tt}-\Div(\TT-\theta\Id)&=f\,,\\
\C^{-1}\TT_t+\xi&=\eps(u_t)\,,\\
\xi&\in N_{K(\theta)}(\TT)\,,\\
\theta_t-\Delta\theta+\theta\Div u_t&=\ip{\xi}{\TT}\,.
\end{aligned}
\end{equation}
Here $\C$ is the elasticity tensor, $\eps(v)=\frac12(\nabla v+\nabla v^\top)$ is the symmetric gradient, and $N_{K(\theta)}$ is the normal cone to the yield set $K(\theta)=\R\Id+r(\theta)C_0$. The set $C_0$ is a compact convex body in the deviatoric stress space, while $r$ is continuous, bounded, nonnegative, and vanishes above a finite critical temperature $\theta_c$. The particular choice $C_0=\overline B_1^{\mathcal S^3_{\rm dev}}(0)$, where $\overline B_1^{\mathcal S^3_{\rm dev}}(0)$ denotes the closed unit ball in $\mathcal S^3_{\rm dev}$, gives the von Mises yield section. Taking
$$
r(\theta)=\{k-\theta\}_+\,,
\qquad
\{s\}_+:=\max\{s,0\}\,,\quad s\in\R\,,
$$
with $k$ playing the role of the critical temperature, gives the explicit yield set 
$$K(\theta)=\{\mathbb S\in\Ssym:\ |\dev\mathbb S|\le\{k-\theta\}_+\}\,.$$
Thus its deviatoric section shrinks linearly for $\theta<k$, collapses to the origin at $\theta=k$, and remains collapsed for $\theta\ge k$. Equivalently, the full yield set contracts onto the hydrostatic axis $\R\Id$. This idealization concerns the loss of shear resistance within a small-strain model; it does not describe a solid--liquid phase transition, latent heat, or a liquid constitutive response.

Low-melting-point metals and alloys provide qualitative motivation for extreme thermal softening. Field's metal, a eutectic bismuth--indium--tin alloy, melts at approximately $62^\circ\mathrm C$ \cite{IdrusSaidi2019}. Experiments document its temperature and strain-rate dependence and strain softening \cite{NguyenDengZhang2020}, without validating the rate-independent law or the particular truncated radius considered here. Pure gallium melts at approximately $29.7646^\circ\mathrm C$ \cite{Strouse2004}, and phase-change materials are used in variable-stiffness robotic systems \cite{GairaEtAl2026}.

The geometry of the degenerating yield set and the mechanical meaning of complete loss of deviatoric strength are discussed in Section~\ref{sec:model} and illustrated in Figure~\ref{fig:thermal-yield-collapse}.

The normal-cone formulation and the variational theory of perfect plasticity originate in the work of Duvaut--Lions, Suquet, and Temam; see \cite{DuvautLions,Suquet1981,Temam1983} and the monographs \cite{HanReddy1999,SimoHughes1998}. Temam's generalized Norton--Hoff approximation constructs the isothermal Prandtl--Reuss evolution and its normality relation \cite{Temam1986}. Anzellotti--Luckhaus \cite{AnzellottiLuckhaus1987} establish a dynamical theory of elasto-perfect plasticity, while Dal Maso--DeSimone--Mora \cite{DalMasoDeSimoneMora2006} develop a quasistatic variational formulation with measure-valued plastic strain. Babadjian--Mifsud \cite{BabadjianMifsud2017} analyse the hyperbolic structure and boundary conditions of a simplified dynamical model. These theories provide relevant reference points for stress--plastic-strain duality and boundary relaxation, although their settings differ from the thermally degenerating system below. In dynamic isothermal plasticity, Kisiel--Che\l{}mi\'nski \cite{KisielChelminski2020} prove existence without safe-load conditions. Here the absence of a safe-load condition is combined with spatial thermal coupling and complete degeneration of the moving yield set.

Several thermoplastic existence theories obtain coercivity from hardening or viscosity. Che\l{}mi\'nski and Racke treat thermoplasticity with linear kinematic hardening \cite{ChelminskiRacke2006}; renormalized formulations for Norton--Hoff thermoviscoplasticity are developed in \cite{ChelminskiOwczarek2016I,ChelminskiOwczarek2016II}. Bartels--Roub\'i\v{c}ek \cite{BartelsRoubicek2008} consider temperature-dependent homothetic yield sets together with Kelvin--Voigt viscosity and a positive-definite rate-dependent term in the plastic evolution. Their yield factor is bounded above and away from zero, and an $L^2$ estimate on the plastic rate supports the limiting constitutive duality argument. 

In \cite{BartelsRoubicek2011}, Bartels--Roub\'i\v{c}ek couple thermal expansion, viscous elasticity, hardening, and rate-independent plasticity with a fixed yield criterion. Related thermodynamic analyses appear in \cite{Roubicek2010,RoubicekPerfect2013}, and vanishing-hardening approximations in \cite{BartelsMielkeRoubicek2012}. Rossi \cite{Rossi2018} develops entropic and weak-energy formulations for thermoviscoplasticity and a joint vanishing-viscosity and vanishing-inertia limit leading to quasistatic perfect plasticity with spatially constant temperature. The limit considered below retains inertia and spatially varying thermal fields while allowing complete collapse of the yield section.

Our approximation retains inertia and the thermal expansion term but introduces neither viscous elasticity nor hardening. Its plastic regularization disappears in the limit. The collapse $r(\theta)=0$ is admissible and no positive lower bound for $r$ is imposed. The associated plastic work controls $r(\theta)|\xi|$, whereas the unweighted multiplier need not be bounded in $L^1$. This distinction determines both the compactness argument and the form of the constitutive information retained by the limit.

Entropy variables provide a second part of the construction. The identity for $\tau=\log\theta$ controls the logarithmic gradient and plastic entropy production even when the heat source is only weakly controlled. Entropy--energy formulations with thermal defects have been developed for nonlinear thermoelasticity in \cite{CieslakMuhaTrifunovic2024} and for Mr\'oz-type inelastic response in \cite{CieslakOwczarekWielgos2026}. In the latter setting the constitutive map has uniform linear growth in stress, which bounds the inelastic rate in $L^2$ and permits a Young-measure description. The degenerating normal-cone constraint considered here has a different coercivity structure: finite measures arise after weighting the multiplier by the yield scale. Related entropy methods in thermodynamics are developed in \cite{FeireislNovotny2009}. The general theories of concentration and generalized Young measures \cite{AlibertBouchitte1997,DalMasoDeSimoneMoraMorini2007} provide a useful context for this distinction. The use of logarithmic entropy also differs from directly identifying the heat equation with an $L^1$ source by truncation methods, as in \cite{BoccardoGallouet1989,Blanchard1993}.

The thermoelastic part of the system belongs to a line of hyperbolic--parabolic
models initiated, in one space dimension, by the global theories of Slemrod,
Dafermos, and Dafermos--Hsiao
\cite{Slemrod1981,Dafermos1982,DafermosHsiao1982}.  Recent results based more
directly on the energy and entropy structures include the regular and
large-time theories of Bies and Cie\'slak
\cite{BiesCieslak2023,BiesCieslak2025} and Winkler's rough-data weak theory
\cite{Winkler2026}.  Multidimensional thermoviscoelastic systems with
$L^1$-temperature data, Kelvin--Voigt rheology, or nonlinear coupling are
treated in
\cite{Roubicek2009,RoubicekNonlinear2013,MielkeRoubicek2020,PawlowZajaczkowski2013}.
These works clarify the energy--entropy compactness available from the thermal
subsystem; the present difficulty is the additional loss of coercivity caused
by complete collapse of the plastic yield section.

The main result constructs global approximation-generated relaxed entropy--energy states and describes the constitutive measures retained in the limit. First, a displacement Galerkin approximation is coupled to a pointwise stress evolution and a full parabolic heat equation. Exact approximate energy and entropy identities yield strong compactness of the logarithmic temperature. Second, a square-root plastic-work estimate yields the rate $\sqrt{\log n/n}$ for the distance to the moving yield set and identifies admissibility of the limiting stress. Third, the support-weighted plastic measures possess lifts concentrated on the exact normality graph, including through complete collapse. Homothetic comparison fields also give a weighted Minty inequality and an anisotropic measure bound. A positive-radius lift organizes the nondegenerate levels without further extraction at its continuity levels.

The approximation framework developed here is also intended to provide the analytical basis for a structure-preserving numerical discretization of the thermoplastic system. A numerical realization of the scheme, together with a computational investigation of the corresponding convergence mechanisms, will be presented in a separate article.

The analysis also identifies the precise scope of the constitutive information retained in the limit. At the macroscopic level, the intrinsic balance and admissibility conditions alone do not imply the Prandtl--Reuss normality relation. Proposition~\ref{prop:classical-consistency} shows that classical solutions belong to the relaxed class, while Proposition~\ref{prop:smooth-nonnormal} provides an explicit smooth state satisfying the intrinsic conditions but violating normality. Thus the gap between the intrinsic relaxed formulation and macroscopic Prandtl--Reuss normality is genuine rather than merely a gap in the compactness argument. Section~\ref{sec:normality} further clarifies the sharpness of this formulation. In the fully collapsed regime $r\equiv0$, an energy--entropy rigidity argument yields approximation-generated states for which the full plastic flow need not be a finite measure or the displacement need not retain its full boundary trace. A simple barycentric argument also shows that exact normality on the lifted fibres does not by itself descend to a macroscopic normality relation.

The model and its thermodynamic identities are specified in Section~\ref{sec:model}. Sections~\ref{sec:state} and \ref{sec:main} give the intrinsic relaxed formulation and the main theorem. Fixed-level solvability, uniform estimates, and the limit equations are proved in Sections~\ref{sec:approximation}--\ref{sec:limit}. Section~\ref{sec:plastic} establishes the weighted and lifted constitutive structure. Its relation to classical normality is described in Section~\ref{sec:normality}. A short geometric appendix records the part of the argument independent of homothety.

\section{Model and thermodynamic structure}
\label{sec:model}
We now specify the mathematical setting for the thermoplastic system \eqref{eq:intro-model}. Let $\Omega\subset\R^3$ be a bounded domain with smooth boundary and let $T>0$. We set
$$
Q=(0,T)\times\Omega\,,\qquad\qbar=[0,T]\times\overline\Omega\,.
$$
We write $\Ssym$ for the space of symmetric $3\times3$ matrices and $\Sdev$ for its trace-free subspace, both endowed with the Frobenius scalar product. For $S\in\Ssym$ and $v:\Omega\to\R^3$, we use the notation
$$
\dev S=S-\frac13(\tr S)\Id\,, \qquad \eps(v)=\frac12(\nabla v+\nabla v^\top)\,,
$$
and tensor divergences are taken row by row. The system \eqref{eq:intro-model} is supplemented with the boundary conditions 
\begin{align}
u&=0 &&\text{on }(0,T)\times\partial\Omega\,,\nonumber\\
\partial_\nu\theta&=0 &&\text{on }(0,T)\times\partial\Omega\,,\label{eq:bc - data}
\end{align}
and initial conditions 
\begin{equation}
u(0)=u_0\,,\qquad u_t(0)=u_1\,,\qquad \mathbb T(0)=\mathbb T_0\,,\qquad \theta(0)=\theta_0\,.
        \label{eq:initial-conditions}
\end{equation}
The precise structural assumptions on the elasticity tensor, the temperature-dependent yield set, and the data are collected below.
\begin{assumption}\label{ass:standing}
The time-independent elasticity tensor $\C(x):\Ssym\to\Ssym$ is measurable, bounded, self-adjoint, and uniformly positive definite. Its inverse satisfies, for constants $c_*,c^*>0$,
\begin{equation}\label{eq:elasticity}
c_*|S|^2\le \C^{-1}(x)S:S\le c^*|S|^2
\qquad(S\in\Ssym),\quad\text{for a.e. }x\in\Omega.
\end{equation}
The set $C_0\subset\Sdev$ is compact and convex, contains zero in its interior, and satisfies
\begin{equation}\label{eq:body}
\overline B_1^{\Sdev}(0)\subset C_0\subset\overline B_{\kappa_0}^{\Sdev}(0),\qquad \kappa_0\ge1.
\end{equation}
The scale $r:[0,\infty)\to[0,\infty)$ is continuous and satisfies
\begin{equation}\label{eq:radius}
0\le r(\theta)\le r_{\max},\qquad r(\theta)=0\quad\text{for }\theta\ge\theta_c,
\end{equation}
for finite constants $r_{\max}>0$ and $\theta_c>0$. The data satisfy
\begin{equation}\label{eq:data}
\begin{gathered}
f\in L^2(Q;\R^3),\quad u_0\in H^1_0(\Omega;\R^3),\quad u_1\in L^2(\Omega;\R^3),\quad \TT_0\in L^2(\Omega;\Ssym),\\
\theta_0>0\ \text{a.e.},\qquad \theta_0,\log\theta_0\in L^1(\Omega),\qquad
\dev\TT_0\in r(\theta_0)C_0\ \text{a.e.}
\end{gathered}
\end{equation}
\end{assumption}
Condition~\eqref{eq:body} fixes the normalization of $C_0$ and $r$. No smoothness or symmetry of $C_0$ is required. In particular, the geometric argument applies to anisotropic convex bodies and nonsmooth yield criteria. Monotonicity of $r$ is not used analytically; a nonincreasing $r$ is the thermal-softening interpretation. The finite-temperature cutoff in \eqref{eq:radius} will provide the global weighted plastic estimate.

For a closed convex set $C$ in a finite-dimensional Euclidean space $E$, our convention is
\begin{equation}\label{eq}
N_C(S)=\{V\in E:\left\langle V,Y-S\right\rangle\le0,\text{ for all },Y\in C\},
\qquad S\in C\,,
\end{equation}
and $N_C(S)=\varnothing$ for $S\notin C$. Define
\begin{equation}\label{eq:yield-set}
K(\theta)=\R\Id+r(\theta)C_0\,,\qquad h_0(V)=\sup_{Y\in C_0}\left\langle Y,V\right\rangle.
\end{equation}
The function $h_0$ is convex, continuous, positively homogeneous, and satisfies $|V|\le h_0(V)\le\kappa_0|V|$. The inclusion $\xi\in N_{K(\theta)}(\TT)$ is equivalent to
\begin{equation}\label{eq:normal-dev}
\xi\in\Sdev,\quad \dev\TT\in r(\theta)C_0,\quad
\ip{\xi}{\TT}=r(\theta)h_0(\xi).
\end{equation}
Indeed, variation along the hydrostatic axis forces $\tr\xi=0$, and maximization over the remaining section gives the support equality. Thus plastic production is nonnegative. If $r(\theta)>0$ and $\dev\TT$ belongs to the interior of $r(\theta)C_0$, then $\xi=0$. At complete collapse, $K(\theta)=\R\Id$ and
\begin{equation}\label{eq:collapsed-normal}
N_{K(\theta)}(a\Id)=\Sdev,\qquad \ip{\xi}{\TT}=0.
\end{equation}
Collapse therefore removes deviatoric resistance while preserving a well-defined normal cone.

\paragraph{Geometry of thermal softening.}
Write $\TT=a\Id+Z$, where $a=\frac{1}{3}\tr\TT$ and $Z=\dev\TT$.  The set $K(\theta)$ is an unbounded cylinder along the hydrostatic axis $\R\Id$: temperature rescales its deviatoric cross-section $r(\theta)C_0$ but does not impose a yield bound on $a$.  Since the thermal stress $-\theta\Id$ is also hydrostatic, $\dev(\TT-\theta\Id)=Z$. The contraction affects only the admissible shear stress. Hydrostatic stress remains unrestricted by the yield condition, and $\C$ remains uniformly positive definite.

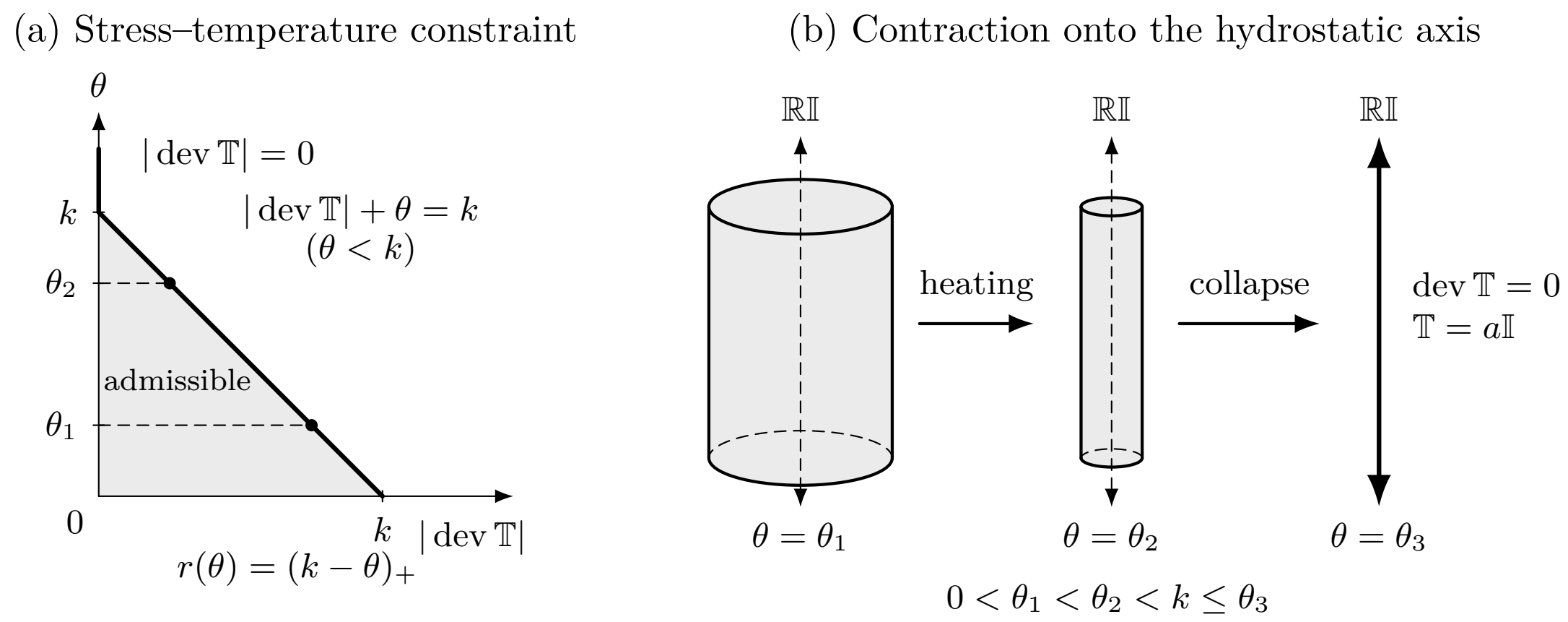
\begin{figure}[!htbp]
\centering
\begin{tikzpicture}[
  x=1cm,y=1cm,>=Latex,
  every node/.style={font=\footnotesize,text=revisiongreen},
  every path/.style={draw=revisiongreen},
  line cap=round,line join=round]
  \node[font=\small] at (2.10,5.05)
    {(a) Stress--temperature constraint};
  \begin{scope}[shift={(0.30,0.80)}]
    \path[fill=revisiongreen!8,draw=none]
      (0,0)--(2.60,0)--(0,2.60)--cycle;
    \draw[->,thin] (0,0)--(3.80,0);
    \node at (3.42,-0.38) {$|\dev\TT|$};
    \draw[->,thin] (0,0)--(0,3.52) node[above] {$\theta$};
    \draw[very thick] (2.60,0)--(0,2.60);
    \draw[very thick] (0,2.60)--(0,3.18);
    \draw[densely dashed,thin] (0,0.65)--(1.95,0.65);
    \draw[densely dashed,thin] (0,1.95)--(0.65,1.95);
    \fill[revisiongreen] (1.95,0.65) circle (1.4pt);
    \fill[revisiongreen] (0.65,1.95) circle (1.4pt);
    \draw[thin] (-0.05,0.65)--(0.05,0.65)
      node[left,xshift=-3pt] {$\theta_1$};
    \draw[thin] (-0.05,1.95)--(0.05,1.95)
      node[left,xshift=-3pt] {$\theta_2$};
    \draw[thin] (-0.05,2.60)--(0.05,2.60)
      node[left,xshift=-3pt] {$k$};
    \draw[thin] (2.60,-0.05)--(2.60,0.05)
      node[below,yshift=-3pt] {$k$};
    \node[below left] at (0,0) {$0$};
    \node[font=\scriptsize] at (0.72,1.07) {admissible};
    \node[align=center] at (2.40,2.43)
      {$|\dev\TT|+\theta=k$\\[-1pt]$(\theta<k)$};
    \node[anchor=west] at (0.25,3.12) {$|\dev\TT|=0$};
  \end{scope}
  \node at (2.12,0.14) {$r(\theta)=(k-\theta)_+$};

  \node[font=\small] at (9.79,5.05)
    {(b) Contraction onto the hydrostatic axis};
  \begin{scope}[shift={(6.73,0)}]
    \path[fill=revisiongreen!8,draw=none]
      (-0.84,1.15)--(-0.84,3.45)
      arc[start angle=180,end angle=0,x radius=0.84,y radius=0.25]
      --(0.84,1.15)
      arc[start angle=0,end angle=-180,x radius=0.84,y radius=0.25]
      --cycle;
    \draw[<->,thin,densely dashed] (0,0.70)--(0,4.10);
    \draw[thick] (-0.84,1.15)--(-0.84,3.45);
    \draw[thick] (0.84,1.15)--(0.84,3.45);
    \draw[thick] (0,3.45) ellipse (0.84 and 0.25);
    \draw[thick] (0.84,1.15)
      arc[start angle=0,end angle=-180,x radius=0.84,y radius=0.25];
    \draw[thin,densely dashed] (0.84,1.15)
      arc[start angle=0,end angle=180,x radius=0.84,y radius=0.25];
    \node[above] at (0,4.10) {$\R\Id$};
    \node at (0,0.43) {$\theta=\theta_1$};
  \end{scope}
  \draw[->,thick] (7.82,2.38)--(8.87,2.38)
    node[midway,above,yshift=2pt] {heating};
  \begin{scope}[shift={(9.58,0)}]
    \path[fill=revisiongreen!8,draw=none]
      (-0.28,1.15)--(-0.28,3.45)
      arc[start angle=180,end angle=0,x radius=0.28,y radius=0.08333]
      --(0.28,1.15)
      arc[start angle=0,end angle=-180,x radius=0.28,y radius=0.08333]
      --cycle;
    \draw[<->,thin,densely dashed] (0,0.70)--(0,4.10);
    \draw[thick] (-0.28,1.15)--(-0.28,3.45);
    \draw[thick] (0.28,1.15)--(0.28,3.45);
    \draw[thick] (0,3.45) ellipse (0.28 and 0.08333);
    \draw[thick] (0.28,1.15)
      arc[start angle=0,end angle=-180,x radius=0.28,y radius=0.08333];
    \draw[thin,densely dashed] (0.28,1.15)
      arc[start angle=0,end angle=180,x radius=0.28,y radius=0.08333];
    \node[above] at (0,4.10) {$\R\Id$};
    \node at (0,0.43) {$\theta=\theta_2$};
  \end{scope}
  \draw[->,thick] (10.20,2.38)--(11.49,2.38)
    node[midway,above,yshift=2pt] {collapse};
  \draw[<->,very thick] (12.03,0.70)--(12.03,4.10);
  \node[above] at (12.03,4.10) {$\R\Id$};
  \node[anchor=west,align=left] at (12.20,2.55)
    {$\dev\TT=0$\\$\TT=a\Id$};
  \node at (12.03,0.43) {$\theta=\theta_3$};
  \node at (9.55,-0.14) {$0<\theta_1<\theta_2<k\le\theta_3$};
\end{tikzpicture}
\caption{Thermally truncated von Mises geometry.
(a) The shaded region and its boundary give the admissible pairs below $k$; above $k$ only the vertical ray $|\dev\TT|=0$ remains.
(b) A schematic section spanned by the hydrostatic direction and two orthonormal deviatoric directions.  The cylinders are unbounded along $\R\Id$; only finite axial portions are shown.  Their lateral boundaries contract with temperature, while the hydrostatic axis survives complete collapse.}
\label{fig:thermal-yield-collapse}
\end{figure}

For the thermally truncated von Mises law, the geometry is especially transparent:
\begin{equation}\label{eq:model-von-mises-green}
K(\theta)=\{a\Id+Z:\ a\in\R,\ Z\in\Sdev,\ |Z|\le\{k-\theta\}_+\}.
\end{equation}
Figure~\ref{fig:thermal-yield-collapse} illustrates the geometry of the truncated von Mises constraint and its contraction with increasing temperature. The particular expression $r(\theta)=\{k-\theta\}_+$ corresponds to normalized variables. More generally, independent reference scales for stress and temperature may be introduced through
$$
r(\theta)=\sigma_{\rm y}^0\left\{1-\frac{\theta}{\theta_c}\right\}_+\,, \qquad
\sigma_{\rm y}^0>0\,,\quad \theta_c>0\,.
$$
The same mechanism is not restricted to circular sections.  After the normalization in \eqref{eq:body}, one may take an ellipsoidal body $C_0=\{Z\in\Sdev:\ip{\mathbb H Z}{Z}\le1\}$, with $\mathbb H$ self-adjoint and positive definite, or the compact convex Tresca body determined by a bound on the difference between the largest and smallest principal stresses. In a principal-stress deviatoric plane the latter has a hexagonal section. In each case the fixed shape is retained as its scale decreases to zero; neither smoothness nor strict convexity of the yield boundary is required.

\paragraph{Mechanical interpretation and limits of the model.}
For a nonincreasing softening law $r$, dissipative heating and contraction of the yield section form a feedback mechanism: plastic work contributes nonnegatively to the heat source, and a resulting temperature increase reduces the range of admissible deviatoric stresses.  This mechanism is already present in thermoviscoplastic models such as \cite{BartelsRoubicek2008}; the distinctive regime studied here permits complete collapse without viscous elasticity or hardening.  Heat diffusion and the expansion term also enter the temperature balance, so heating need not be monotone along a solution.

The loss of coercivity has a direct mechanical interpretation.  At the smooth level, \eqref{eq:normal-dev} gives, for every $\delta>0$,
\begin{equation}\label{eq:model-work-coercivity-green}
\ip{\xi}{\TT}=r(\theta)h_0(\xi)\ge r(\theta)|\xi|\,, \qquad |\xi|\le\delta^{-1}\ip{\xi}{\TT}
\quad\text{where }r(\theta)\ge\delta\,.
\end{equation}
There is no uniform estimate of this form as $r(\theta)\downarrow0$. At complete collapse, the normal cone admits every deviatoric flow direction while the classical plastic work is zero; the actual rate must still satisfy the coupled mechanical equations.  The model thus allows, at the constitutive level, irreversible shear with vanishing plastic resistance. The resulting loss of plastic coercivity leads to the weighted compactness and lifted description used below.

This reduction of shear capacity motivates the study of heated load-bearing components, such as walls or thin material layers. The present small-strain model isolates the constitutive mechanism and does not analyse buckling, folding, fracture, or shear-band formation; those phenomena require further structural modelling, including finite-deformation or shell theories where appropriate. The collapse in Figure~\ref{fig:thermal-yield-collapse} concerns the admissible deviatoric stress section.

For a smooth positive-temperature solution, testing momentum by $u_t$ and the constitutive equation by $\TT$ gives
\begin{equation}\label{eq:mechanical-formal}
\frac{\mathrm d}{\mathrm dt}\left(\frac12\norm{u_t}_{L^2(\Omega)}^2+
\frac12\int_\Omega\ip{\C^{-1}\TT}{\TT}\dd x\right)
=\int_\Omega \ip{f}{u_t}\dd x+\int_\Omega\theta\Div u_t\dd x-\int_\Omega\ip{\xi}{\TT}\dd x\,.
\end{equation}
Integration of the heat equation cancels the last two terms, so
\begin{equation}
\begin{aligned}
\label{eq:energy-formal}
\frac{\mathrm d}{\mathrm dt}E(t)&=\int_\Omega \ip{f}{u_t}\dd x\,,\\
E(t)&=\frac12\norm{u_t(t)}_{L^2(\Omega)}^2+\frac12\int_\Omega\ip{\C^{-1}\TT(t)}{\TT(t)}\dd x+\int_\Omega\theta(t)\dd x\,.
\end{aligned}
\end{equation}
Dividing the heat equation by $\theta$ yields, with $\tau=\log\theta$,
\begin{equation}\label{eq:entropy-formal}
\partial_t(\tau+\Div u)-\Delta\tau
=\frac{\ip{\xi}{\TT}}{\theta}+|\nabla\tau|^2.
\end{equation}
The entropy density is $\tau+\Div u$ and its flux is $-\nabla\tau$. Equation~\eqref{eq:normal-dev} gives nonnegative plastic entropy production. Since $\int_\Omega\Div u=0$, subtraction of the integrated entropy identity from \eqref{eq:energy-formal} gives
\begin{equation}\label{eq:free-energy-formal}
\begin{aligned}
E(t)&-\int_\Omega\tau(t)\dd x+
\int_0^t\!\int_\Omega\left(\frac{\ip{\xi}{\TT}}{\theta}+|\nabla\tau|^2\right)\dd x\dd s
=E(0)\\
&\qquad-\int_\Omega\log\theta_0\dd x+\int_0^t\!\int_\Omega \ip{f}{u_t}\dd x\dd s\,.
\end{aligned}
\end{equation}
The approximation preserves these identities. The scale $r$ enters through the instantaneous dissipation potential $r(\theta)h_0$, so its differentiation is not required.

\section{Relaxed entropy--energy states}
\label{sec:state}
For $s>3$ we set $Y_s=(W^{1,s}(\Omega))^*$. The space $\M(X;E)$ denotes the finite $E$-valued Radon measures on $X$, and $\M^+(X)$ denotes the cone of finite nonnegative Radon measures on
$X$. We write $BV([0,T];Y_s)$ for the space of $Y_s$-valued functions of bounded variation in time; equivalently, $v\in BV([0,T];Y_s)$ if its distributional time derivative $D_tv$ is a finite $Y_s$-valued Radon measure on $(0,T)$.

\noindent
Since $\M(\overline\Omega)=\left(C(\overline\Omega)\right)^*$, the notation $L^\infty_{w*}(0,T;\M(\overline\Omega))$ denotes the space of weak-star measurable, essentially bounded maps
$t\mapsto\mu_t\in\M(\overline\Omega)$. Here weak-star measurability means that $t\longmapsto\langle\mu_t,\varphi\rangle$ is measurable for every $\varphi\in C(\overline\Omega)$, while essential boundedness means
$$
\operatorname*{ess\,sup}_{t\in(0,T)}
\|\mu_t\|_{\M(\overline\Omega)}<\infty\,.
$$
As usual, such maps are identified up to equality for almost every $t\in(0,T)$. For nonnegative measures, $\|\mu_t\|_{\M(\overline\Omega)}=\mu_t(\overline\Omega)$. For a Banach space $X$, we denote by $C_w([0,T];X)$ the space of maps from $[0,T]$ into $X$ that are continuous with respect to the weak topology of $X$; equivalently, $v\in C_w([0,T];X)$ if $t\longmapsto\langle x^*,v(t)\rangle$ is continuous on $[0,T]$ for every $x^*\in X^*$.

\noindent
We use
\[
\begin{aligned}
\Hdiv&=\{v\in L^2(\Omega;\R^3):\Div v\in L^2(\Omega)\},\\
\Hzdiv&=\{v\in\Hdiv:\ip{v}{\nu}=0\text{ in }H^{-1/2}(\partial\Omega)\}.
\end{aligned}
\]
The normal trace is continuous from $\Hdiv$ into $H^{-1/2}(\partial\Omega)$, and $\Hzdiv$ is the closure of $C_c^\infty(\Omega;\R^3)$, and hence also of
$H^1_0(\Omega;\R^3)$, in $\Hdiv$. Consequently, weak $\Hdiv$ limits of
sequences from $H^1_0(\Omega;\R^3)$ retain the zero normal trace, although
the full Dirichlet trace need not survive.

\noindent
If $F:X\to Y$ is a measurable map and $\mu$ is a finite Radon measure on $X$, we denote by $F_\#\mu$ the push-forward of $\mu$ under $F$, defined by
\[
(F_\#\mu)(B)=\mu(F^{-1}(B))
\]
for every Borel set $B\subset Y$. Equivalently,
\[
\int_Y \psi(y)\,\dd(F_\#\mu)(y)
=
\int_X \psi(F(x))\,\dd\mu(x)
\]
for every bounded Borel function $\psi:Y\to\mathbb R$.

\noindent
Set $\tau_0=\log\theta_0$ and
\begin{equation}\label{eq:E0}
E_0=\frac12\norm{u_1}_2^2+\frac12\int_\Omega\ip{\C^{-1}\TT_0}{\TT_0}\dd x+\int_\Omega\theta_0\dd x.
\end{equation}

\begin{definition}[Relaxed entropy--energy state]\label{def:state}
Under Assumption~\ref{ass:standing}, a relaxed entropy--energy state is a tuple
\[
(u,\TT,\tau,\Theta,\Xi,\mu_{\rm pl},\mu_{\rm th})
\]
with the following properties.
\begin{enumerate}
\item Its regularity is
\begin{equation}\label{eq:state-regularity}
\begin{aligned}
u&\in W^{1,\infty}(0,T;L^2(\Omega;\R^3))\cap L^\infty(0,T;\Hzdiv)\,,\\
u_t&\in C_w([0,T];L^2(\Omega;\R^3))\,,\qquad
\TT\in L^\infty(0,T;L^2(\Omega;\Ssym))\,,\\
\tau&\in L^\infty(0,T;L^1(\Omega))\cap L^2(0,T;H^1(\Omega))\cap BV([0,T];Y_s)\,,\\
e^\tau&\in L^\infty(0,T;L^1(\Omega))\,,\qquad
\Theta\in L^\infty_{w*}(0,T;\M^+(\overline\Omega))\,,\\
\Xi&\in H^{-1}(Q;\Sdev)\,,\qquad \mu_{\rm pl},\mu_{\rm th}\in\M^+(\qbar)\,.
\end{aligned}
\end{equation}
for every finite $s>3$. There exists $\lambda\in L^\infty_{w*}(0,T;\M^+(\overline\Omega))$ such that
\begin{equation}\label{eq:thermal-split}
\Theta(t)=e^{\tau(t)}\dd x+\lambda_t\quad\text{for a.e. }t\in(0,T).
\end{equation}

\item The initial displacement and velocity are $u(0)=u_0$ strongly in $L^2(\Omega;\R^3)$ and $u_t(0)=u_1$ weakly in $L^2(\Omega;\R^3)$. For every $\varphi\in C_c^\infty([0,T)\times\Omega;\R^3)$,
\begin{equation}\label{eq:state-momentum}
\begin{split}
-\int_Q \ip{u_t}{\varphi_t}\dd x\dd t+\int_Q \ip{\TT}{\eps(\varphi)}\dd x\dd t
&\\[-1mm]\quad-\int_0^T \dual{\Theta(t)}{\Div\varphi(t)}\dd t &=
\int_Q \ip{f}{\varphi}\dd x\dd t +\int_\Omega \ip{u_1}{\varphi(0)}\dd x\, .
\end{split}
\end{equation}

\item Admissibility and the distributional constitutive balance hold:
\begin{equation}\label{eq:state-constitutive}
\dev\TT\in r(e^\tau)C_0\quad\text{a.e. in }Q,\qquad
\partial_t(\C^{-1}\TT)+\Xi=\eps(u_t)\quad\text{in }\mathcal D'(Q;\Ssym).
\end{equation}
The initial-time extension $\overline\Xi$ is defined, for symmetric $\Phi\in C_c^\infty([0,T)\times\Omega;\Ssym)$, by
\begin{equation}\label{eq:initial-constitutive}
\dual{\overline\Xi}{\Phi} = -\int_Q \ip{u_t}{\Div\Phi}\dd x\dd t +\int_Q \ip{\C^{-1}\TT}{\Phi_t}\dd x\dd t +\int_\Omega \ip{\C^{-1}\TT_0}{\Phi(0)}\dd x.
\end{equation}
This extension is trace-free, and
\begin{equation}\label{eq:state-volume}
\Div u-\tr(\C^{-1}\TT)=\Div u_0-\tr(\C^{-1}\TT_0)\quad\text{a.e. in }Q.
\end{equation}

\item For every $\phi\in C_c^\infty([0,T)\times\overline\Omega)$,
\begin{equation}\label{eq:state-entropy}
\begin{aligned}
-\int_Q (\tau+\Div u)\phi_t\,\dd x\dd t+\int_Q \ip{\nabla\tau}{\nabla\phi}\,\dd x\dd t &=
\int_{\qbar}\phi\,\dd(\mu_{\rm pl}+\mu_{\rm th})\\
&\quad
+\int_\Omega (\tau_0+\Div u_0)\phi(0)\,\dd x
\end{aligned}
\end{equation}
and $\mu_{\rm th}\ge |\nabla\tau|^2\dd x\dd t$ as measures on $\qbar$, extending the density by zero outside $Q$.

\item With
\begin{equation}\label{eq:state-energy}
\mathcal E(t)=\frac12\norm{u_t(t)}_2^2+
\frac12\int_\Omega\ip{\C^{-1}\TT(t)}{\TT(t)}\dd x+\Theta(t)(\overline\Omega),
\end{equation}
the inequalities
\begin{equation}\label{eq:state-energy-ineq}
\mathcal E(t)\le E_0+\int_0^t\!\int_\Omega \ip{f}{u_t}\dd x\dd \tau
\end{equation}
and
\begin{equation}\label{eq:state-combined}
\begin{aligned}
\mathcal E(t)-\int_\Omega\tau(t)\dd x
+\mu_{\rm pl}([0,t]\times\overline\Omega)&+\mu_{\rm th}([0,t]\times\overline\Omega)
\le E_0\\
&\qquad-\int_\Omega\tau_0\dd x+\int_0^t\!\int_\Omega \ip{f}{u_t}\dd x\dd \tau
\end{aligned}
\end{equation}
hold for almost every $t\in(0,T)$.
\end{enumerate}
\end{definition}
The spatial boundary condition for $u$ is its zero normal trace in $\Hzdiv$. Thermal insulation is encoded by the entropy test class reaching $\partial\Omega$ in \eqref{eq:state-entropy}. This formulation accommodates boundary concentrations of entropy production. The initial logarithmic temperature and thermal energy enter through \eqref{eq:state-entropy} and $E_0$, respectively, and production atoms at $t=0$ are retained. The stress datum is encoded in \eqref{eq:initial-constitutive} and \eqref{eq:state-volume}. No separate strong temperature trace or full $L^2$ stress trace is prescribed. The thermal excess $\lambda_t$ is not necessarily singular with respect to Lebesgue measure. It may contain an absolutely continuous part.

\begin{remark}[Integrated entropy and the volumetric stress trace]\label{rem:entropy-volume}
Testing \eqref{eq:state-entropy} by spatial constants and using $\int_\Omega\Div u(t)=0$ gives
\begin{equation}\label{eq:state-integrated-entropy}
\int_\Omega\tau(t)\dd x=\int_\Omega\tau_0\dd x+
(\mu_{\rm pl}+\mu_{\rm th})([0,t]\times\overline\Omega)
\quad\text{for a.e. }t<T.
\end{equation}
Consequently, under the entropy identity, \eqref{eq:state-combined} and \eqref{eq:state-energy-ineq} are equivalent. The combined form corresponds to the estimate used in the approximation and imposes no additional constitutive selection. Moreover \eqref{eq:state-volume}, continuity of $u$ in $L^2$, and the uniform $L^2$ bound on its divergence give a representative
\[
\tr(\C^{-1}\TT)\in C_w([0,T];L^2(\Omega)),\qquad
\tr(\C^{-1}\TT)(0)=\tr(\C^{-1}\TT_0).
\]
Indeed the volume identity first gives continuity in $H^{-1}(\Omega)$, and density upgrades it to weak $L^2$ continuity. This conclusion concerns only the volumetric component and supplies no additional trace of the full stress.
\end{remark}

\begin{remark}[Scope of the intrinsic formulation]\label{rem:scope}
Definition~\ref{def:state} does not impose the macroscopic Prandtl--Reuss normality relation $
\Xi\in N_{K(e^\tau)}(\TT)$. This is why we use the term \emph{relaxed entropy--energy state} rather
than a weak Prandtl--Reuss solution. Every sufficiently regular classical solution induces such a state, as shown in Proposition~\ref{prop:classical-consistency}, but the converse is false: the intrinsic formulation admits even smooth states with nonzero constitutive flow while the stress lies strictly inside a nondegenerate yield set, see Proposition~\ref{prop:smooth-nonnormal}. Theorem~\ref{thm:main} constructs a distinguished subclass arising from
the prescribed approximation, together with additional weighted and
lifted constitutive measures. This approximation-generated subclass is
specified in Definition~\ref{def:generated}. The relation between the
lifted normality retained by these measures and macroscopic
Prandtl--Reuss normality is discussed in Section~\ref{sec:normality}.
\end{remark}

\begin{proposition}[Consistency with classical solutions]
\label{prop:classical-consistency}
Let $(u,\TT,\theta,\xi)$ be a sufficiently regular solution of \eqref{eq:intro-model} with the boundary and initial conditions \eqref{eq:bc - data} and \eqref{eq:initial-conditions}, with $\theta>0$, such that
\[
\TT\in K(\theta),\qquad \xi\in N_{K(\theta)}(\TT)
\quad\text{a.e. in }Q\,.
\]
Then it induces a relaxed entropy--energy state by setting
\[
\tau=\log\theta\,,\qquad \Theta(t)=\theta(t)\dd x\,,\qquad \Xi=\xi\,,
\qquad
\mu_{\rm pl}=\frac{\ip{\xi}{\TT}}{\theta}\dd x\dd t,\qquad
\mu_{\rm th}=|\nabla\log\theta|^2\dd x\dd t\,.
\]
All concentration defects vanish, the inequalities \eqref{eq:state-energy-ineq}--\eqref{eq:state-combined} are equalities, and the constitutive lift is monokinetic. Choose any measurable unit deviator field $\omega$ satisfying $\omega=\xi/|\xi|$ on $\{\xi\ne0\}$, and define
\[
\Gamma_{\rm pl}
:=(\operatorname{id},r(\theta),\dev\TT,\omega)_\#\mu_{\rm pl}\,.
\]
The marginal of $\Gamma_{\rm pl}$ is $\mu_{\rm pl}$, and $\Gamma_{\rm pl}$ is concentrated on the set of fibre states $(\rho,Z,\omega)$ satisfying
\[
Z\in\rho C_0\,,\qquad \ip{Z}{\omega}=\rho h_0(\omega)\,,\qquad |\omega|=1\,.
\]
More precisely, wherever $\xi\ne0$,\( (\rho,Z,\omega)=\left(r(\theta),\dev\TT,\frac{\xi}{|\xi|}\right)\).
\end{proposition}
\begin{proof}
The weak momentum and constitutive balances follow directly from the
classical equations by integration by parts. Since $\theta>0$, division of
the heat equation by $\theta$ gives
$$
\partial_t(\log\theta+\Div u)-\Delta\log\theta=\frac{\ip{\xi}{\TT}}{\theta}+|\nabla\log\theta|^2 $$
and hence the entropy identity with
$$
\mu_{\rm pl}=
\frac{\ip{\xi}{\TT}}{\theta}\dd x\dd t\,,\qquad\mu_{\rm th}=|\nabla\log\theta|^2\dd x\dd t\,.
$$
The classical total-energy identity yields equality in \eqref{eq:state-energy-ineq}, and therefore also in \eqref{eq:state-combined}. Moreover, $\Theta(t)=\theta(t)\dd x=e^{\tau(t)}\dd x$, so the thermal excess vanishes. By the normality relation and \eqref{eq:normal-dev}, on ${\xi\ne0}$,
with $\omega=\xi/|\xi|$,
$$
\ip{\xi}{\TT}=\ip{\xi}{\dev\TT}=r(\theta)h_0(\omega)|\xi|\ge0\,.
$$
Thus $\mu_{\rm pl}$ is a nonnegative measure. Furthermore,
$$
\dev\TT\in r(\theta)C_0\,,\qquad \ip{\dev\TT}{\omega} = r(\theta)h_0(\omega)
$$
wherever $\xi\ne0$, so the map $(t,x)\longmapsto\bigl(t,x,r(\theta),\dev\TT,\omega\bigr)$ takes $\mu_{\rm pl}$-almost every point into the exact normality graph. Consequently, $ \Gamma_{\rm pl}=
(\operatorname{id},r(\theta),\dev\TT,\omega)_\#\mu_{\rm pl}$ is concentrated on that graph. Its projection onto the spacetime variables is $\mu_{\rm pl}$, since the composition of the projection with the above map is the identity.

\noindent
Finally, $\mu_{\rm pl}$ vanishes on ${\xi=0}$, so the arbitrary choice
of the unit field $\omega$ there does not affect $\Gamma_{\rm pl}$.
Notice that $\mu_{\rm pl}$ may also vanish where $\xi\ne0$ in the
collapsed regime $r(\theta)=0$; such zero-work plastic flow is therefore
not recorded by this entropy-weighted lift.
\end{proof}

\section{Main result}
\label{sec:main}
We first fix the notation for the constitutive conclusions. Let $P\in\M(\qbar;\Sdev)$ and write its polar decomposition as
\[
P=\zeta\,|P|\,,\qquad\zeta=\frac{\dd P}{\dd|P|}\,,\qquad
|\zeta|=1\quad |P|\text{-a.e.}
\]
We define its $C_0$-anisotropic variation by
\begin{equation}\label{eq:anisotropic-variation}
|P|_{C_0}=h_0(\zeta)\,|P|\,.
\end{equation}
This is a nonnegative finite measure; symmetry of $C_0$ is not required.

\noindent
Let $\mathcal U=\{\omega\in\Sdev:|\omega|=1\}$ and introduce the compact fibre and the normality graph
\begin{equation}\label{eq:global-fibre}
\mathcal F =[0,r_{\max}]\times\overline B_{\kappa_0r_{\max}}^{\Sdev}(0)\times\mathcal U\,,
\qquad \mathcal G=\left\{(\rho,Z,\omega)\in\mathcal F:Z\in\rho C_0,\ \ip{Z}{\omega}=\rho h_0(\omega) \right\}\,.
\end{equation}
We denote by $\pi:\qbar\times\mathcal F\to\qbar$ the canonical projection onto the spacetime variables. For a nonempty closed convex set $C\subset\Sdev$, we denote by $\Pi_C:\Sdev\to C$ the metric projection onto $C$, uniquely characterized by
\[
\Pi_C(Z) =\operatorname*{arg\,min}_{Y\in C}|Z-Y|\,.
\]
For a measurable set $A$, we denote by $\mathbf 1_A$ its indicator function. We are now in a position to state the main result of the paper.
\begin{theorem}[Global construction and constitutive measure structure]\label{thm:main}
Under Assumption~\ref{ass:standing}, there exists a relaxed entropy--energy state in the sense of Definition~\ref{def:state} and a common subsequence of the approximation constructed in Section~\ref{sec:approximation} for which the following properties hold.
\begin{enumerate}
\item Writing $\tau^n=\log\theta^n$, $r_n=r(\theta^n)$, and $r=r(e^\tau)$, one has
\begin{equation}\label{eq:main-tau}
\begin{aligned}
\tau^n&\longrightarrow\tau &&\text{strongly in }L^2(0,T;L^q(\Omega))\,,\quad 1\le q<6\,,\\
\tau^n&\longrightarrow\tau &&\text{strongly in }L^p(Q),\quad 1\le p<8/3\,,\\
\theta^n&\longrightarrow e^\tau &&\text{a.e. in }Q\,,\\
r_n&\longrightarrow r &&\text{strongly in }L^p(Q),\quad 1\le p<\infty\,.
\end{aligned}
\end{equation}
Moreover $\theta^n\dd x\weakstar\Theta$ in $L^\infty_{w*}(0,T;\M(\overline\Omega))$.

\item For $Z_n=\Pi_{r_nC_0}(\dev\TT^n)$ and $d_n=|\dev\TT^n-Z_n|$, each fixed $1\le p<\infty$ satisfies
\begin{equation}\label{eq:main-rate}
\norm{d_n}_{L^p(Q)}\le C_p\sqrt{\frac{\log n}{n}}
\quad\text{for all sufficiently large }n\,.
\end{equation}
The same bound holds for $\widehat\TT^n-\TT^n$, where $\widehat\TT^n=\frac13\tr\TT^n\Id+Z_n$. In addition,
\begin{equation}\label{eq:main-admissibility}
\begin{gathered}
\dev\TT\in rC_0\quad\text{a.e.}\,,\qquad \dev\TT\in L^\infty(Q;\Sdev)\,,\\
\mathbf 1_{\{r(e^\tau)=0\}}\dev\TT^n \longrightarrow0 \quad\text{in every finite }L^p(Q)\,.
\end{gathered}
\end{equation}

\item The full approximate plastic flow converges weakly in $H^{-1}(Q;\Sdev)$ to $\Xi$. The finite measures
\begin{equation}\label{eq:main-weighted}
r_n\xi^n\dd x\dd t\weakstar\Lambda\,,\quad r_n|\xi^n|\dd x\dd t\weakstar\nu\,,\quad (Z_n:\xi^n)\dd x\dd t\weakstar\chi
\end{equation}
on $\qbar$ satisfy
\begin{equation}\label{eq:main-domination}
|\Lambda|\le\nu\le\chi\le\theta_c\mu_{\rm pl}\,,\qquad |\Lambda|_{C_0}\le\chi\,.
\end{equation}
For nonnegative $\varphi\in C(\qbar)$ and continuous fields $\bar Z,\eta$ with $\bar Z\in C_0$ and $\eta\in N_{C_0}(\bar Z)$ pointwise, the weighted Minty inequality holds:
\begin{equation}\label{eq:main-minty}
\int_{\qbar}\varphi\,\dd\chi -\int_{\qbar}\varphi\,\ip{\bar Z}{\dd\Lambda} -\int_Q \varphi\,\ip{\eta}{\dev\TT-r\bar Z}\dd x\dd t \ge0\,.
\end{equation}

\item There exist positive measures $\Gamma_{\rm pl}$ and $\widetilde\Gamma$ on $\qbar\times\mathcal F$, concentrated on $\qbar\times\mathcal G$, such that
\begin{equation}\label{eq:main-lifts}
\begin{gathered}
\pi_\#\Gamma_{\rm pl}=\mu_{\rm pl}\,,\qquad \pi_\#\widetilde\Gamma=\nu\,,\\
\Lambda=\pi_\#(\omega\widetilde\Gamma)\,,\qquad \chi=\pi_\#(h_0(\omega)\widetilde\Gamma)\,,\\
h_0(\omega)\widetilde\Gamma\le\theta_c\Gamma_{\rm pl}\,.
\end{gathered}
\end{equation}
Thus $\ip{Z}{\omega}=\rho h_0(\omega)=\max_{S\in\rho C_0}\ip{S}{\omega}$ almost everywhere for both lifts. The plastic entropy measure also splits as $\mu_{\rm pl}=\mu_{\rm proj}+\mu_{\rm over}$ with nonnegative components obtained from projected and overstress production.

\item Set $X_+=(\qbar\times\mathcal F)\cap\{\rho>0\}$. The measure
\begin{equation}\label{eq:main-positive-lift}
\Gamma^+=\rho^{-1}\widetilde\Gamma\restriction_{X_+}
\end{equation}
is positive, locally finite on $X_+$, and concentrated on the normality graph. For every $\delta\in(0,r_{\max})$ such that $\widetilde\Gamma(\{\rho=\delta\})=0$, the same subsequence, without further extraction, satisfies
\begin{equation}\label{eq:main-local}
\begin{aligned}
\mathbf1_{\{r_n\ge\delta\}}\xi^n\dd x\dd t&\weakstar P_\delta\,,\\
\mathbf1_{\{r_n\ge\delta\}}(Z_n:\xi^n)\dd x\dd t&\weakstar\mu_\delta\,.
\end{aligned}
\end{equation}
With $\Gamma_\delta=\Gamma^+\restriction_{\{\rho\ge\delta\}}$, regarded as a finite measure on $\qbar\times\mathcal F$, one has
\begin{equation}\label{eq:main-local-lift}
\begin{gathered}
P_\delta=\pi_\#(\omega\Gamma_\delta)\,,\qquad \mu_\delta=\pi_\#(\rho h_0(\omega)\Gamma_\delta)\,,\\
\delta|P_\delta|_{C_0}\le\mu_\delta\le\chi\,.
\end{gathered}
\end{equation}
There are at most countably many exceptional positive levels. For any prescribed countable set $D\subset(0,r_{\max})$, a further common extraction yields graph-supported lifts $\Gamma_\delta$ and the same convergence and bounds for every $\delta\in D$, including exceptional levels.
\end{enumerate}
\end{theorem}

\noindent
At an exceptional positive level $\delta$, the extracted lift need not coincide with the canonical restriction $\Gamma^+\restriction_{\{\rho\ge\delta\}}$, because the limiting support-weighted measure may carry mass on the level set ${\rho=\delta}$. Moreover, the measures in \eqref{eq:main-local} are obtained from the approximate regions ${r_n\ge\delta}$. In general, they cannot be identified with restrictions to the measurable limiting set ${r(e^\tau)\ge\delta}$. A stronger identification is available on fixed open regions where the approximate radii remain uniformly positive. If $U\subset Q$ and $r_n\ge\delta>0$ a.e. on $U$ for all sufficiently large $n$, then the full distributional flow $\Xi\restriction_U$ is a finite measure and coincides with every selected level-restricted limit $P_\eta\restriction_U$ for $0<\eta\le\delta$, see Proposition~\ref{prop:fixed-region-identification}.

\noindent
All constants in the uniform estimates depend only on the data, the domain,
$T$, and the standing structural constants, with the indicated dependence
on $p$ or $\delta$.\\

\noindent
The proof is organized in four steps. We first construct the approximate
solutions and derive the exact energy and logarithmic entropy identities.
The resulting uniform estimates provide the quantitative overstress bound
and the compactness of the logarithmic temperature. These estimates are
then used to pass to the limit in the intrinsic balance laws and to prove
items~\textup{(i)} and \textup{(ii)}. The remaining constitutive information
is recovered through support-weighted compactness and a weighted Minty
argument, yielding item~\textup{(iii)}, followed by the construction of the
global and positive-radius lifts in items~\textup{(iv)} and \textup{(v)}.
The details are carried out in the next four sections.

\section{Approximation and fixed-level solvability}
\label{sec:approximation}

\subsection{Projection geometry and the Norton--Hoff operator}

For $\theta\ge0$ and $S\in\Ssym$, put
\begin{equation}\label{eq:projection-definition}
Z(\theta,S)=\Pi_{r(\theta)C_0}(\dev S)\,,\qquad
R(\theta,S)=\dev S-Z(\theta,S)\,,\qquad d(\theta,S)=|R(\theta,S)|\,.
\end{equation}
For integers $n\ge1$ define
\begin{equation}\label{eq:operator}
a_n(s)=s(1+s^2)^{(n-1)/2},\qquad
G_n(\theta,S)=(1+d^2(\theta,S))^{(n-1)/2}R(\theta,S)\,.
\end{equation}
This is a distance penalization of the moving constraint. It is used only in the construction; $n$ tends to infinity in the final result.

\begin{lemma}[Projection and operator properties]\label{lem:geometry}
The map $G_n$ is continuous in $(\theta,S)$, monotone in $S$ for fixed $\theta$, and locally Lipschitz in $S$, with a Lipschitz constant on bounded stress sets independent of $\theta$. It takes values in $\Sdev$. If $d>0$ and $\omega=R/d$, then
\begin{equation}\label{eq:projection-normality}
\omega\in N_{r(\theta)C_0}(Z)\,,\qquad \ip{Z}{\omega}=r(\theta)h_0(\omega)\,.
\end{equation}
Consequently, with $\xi=G_n(\theta,S)$ and $q=\ip{\xi}{S}$,
\begin{equation}\label{eq:work-identity}
\begin{aligned}
\xi&\in N_{r(\theta)C_0}(Z)\,,\qquad \ip{Z}{\xi}=r(\theta)h_0(\xi)\,,\\
q&=|\xi|\bigl(r(\theta)h_0(\omega)+d\bigr)\ge d^2(1+d^2)^{(n-1)/2}\ge0\,.
\end{aligned}
\end{equation}
Expressions weighted by $|\xi|$ are assigned zero when $d=0$.
\end{lemma}

\begin{proof}
The characterization of metric projection is $\ip{\dev S-Z}{Y-Z}\le0$ for $Y\in r(\theta)C_0$. Dividing by $d$ gives normality and maximization of $\ip{Y}{\omega}$, proving \eqref{eq:projection-normality}. Since normal cones are cones, the same relation holds for $\xi=a_n(d)\omega$. Decomposing $\dev S=Z+d\omega$ proves \eqref{eq:work-identity}.

\noindent
For fixed $\theta$, the function
\[
S\longmapsto A_n(d(\theta,S)),\qquad A_n(s)=\frac{(1+s^2)^{(n+1)/2}-1}{n+1},
\]
is convex: distance to a convex set is convex, and $A_n$ is convex and nondecreasing on $[0,\infty)$. Its gradient is $G_n$, including at $d=0$, by differentiability of the squared distance; see also the monotone-operator framework in \cite{Roubicek2005}. Thus $G_n$ is monotone. The residual $Y-\Pi_CY$ is nonexpansive for a closed convex set. Because $0\in r(\theta)C_0$, its norm is at most $|Y|$. Composition with the radial locally Lipschitz map $V\mapsto(1+|V|^2)^{(n-1)/2}V$ proves the uniform-on-bounded-stress Lipschitz assertion.

\noindent
Finally, $r(\theta_j)C_0$ converges in Hausdorff distance to $r(\theta)C_0$ when $\theta_j\to\theta$. If also $S_j\to S$, compactness of the projected points and passage to the projection inequality identify every cluster point with the unique projection of $\dev S$. This proves joint continuity, including at $r(\theta)=0$.
\end{proof}

\subsection{Data and the semidiscrete system}
Let \((w_j)_{j\ge1}\) be an orthonormal basis of \(L^2(\Omega;\mathbb R^3)\), consisting of eigenfunctions of the vector Dirichlet Laplacian, and set $W_n=\operatorname{span} \{w_1,\ldots,w_n\}$. Let \(P_n:L^2(\Omega;\mathbb R^3)\to W_n\) denote the corresponding \(L^2(\Omega;\mathbb R^3)\)-orthogonal projection. Then $P_nu_0\to u_0\quad\text{in }H^1_0(\Omega;\mathbb R^3)$ and $P_nu_1\to u_1\quad\text{in }L^2(\Omega;\mathbb R^3)$. For every smooth vector field compactly supported in \(\Omega\), the spectral projections converge in every fixed Sobolev norm, uniformly on compact time intervals when applied to time-dependent test fields, also after time differentiation. In particular, they converge in \(W^{1,\infty}(\Omega;\mathbb R^3)\).

\noindent
Choose smooth Neumann-compatible temperatures $\theta_{0,n}$ satisfying
\begin{equation}\label{eq:initial-temp-approx}
0<c_n\le\theta_{0,n}\le C_n\,,\qquad \theta_{0,n}\to\theta_0\,,\quad\log\theta_{0,n}\to\log\theta_0
\quad\text{in }L^1(\Omega)\text{ and a.e.}
\end{equation}
Such data are obtained by truncating $\theta_0$ to $[1/m,m]$ and applying the Neumann heat semigroup for a positive time tending to zero. For fixed $m$, strong $L^1$ continuity of the semigroup and Lipschitz continuity of the logarithm on $[1/m,m]$ give simultaneous approximation of the temperature and its logarithm. The maximum principle preserves these bounds and the smoothing gives Neumann compatibility. A diagonal choice, followed if necessary by a subsequence, proves \eqref{eq:initial-temp-approx}.

\noindent
Bounded admissible stresses can be chosen explicitly. If $b_n$ is the truncation of $\tr\TT_0$ to $[-n,n]$, let
\begin{equation}\label{eq:initial-stress-approx}
\TT_0^n=\frac13b_n\Id+\Pi_{r(\theta_{0,n})C_0}(\dev\TT_0)\,.
\end{equation}
Then $\TT_0^n\in L^\infty(\Omega;\Ssym)$, $\dev\TT_0^n\in r(\theta_{0,n})C_0$, and $\TT_0^n\to\TT_0$ in $L^2(\Omega;\Ssym)$. Indeed, initial admissibility bounds $\dev\TT_0$ by $\kappa_0r_{\max}$, continuity of projection and dominated convergence treat the deviatoric part and scalar truncation treats the trace. Set $u_0^n=P_nu_0$ and $u_1^n=P_nu_1$.

\noindent
We seek $u^n:[0,T]\to W_n$, a stress $\TT^n:Q\rightarrow\Ssym$, and a temperature $\theta^n:Q\rightarrow (0,\infty)$ satisfying
\begin{subequations}\label{eq:approx}
\begin{align}
\int_\Omega \ip{u_{tt}^n}{w}\dd x +\int_\Omega \ip{\TT^n}{\eps(w)}\dd x -\int_\Omega \theta^n \Div w\dd x &= \int_\Omega \ip{f}{w}\dd x\,,\label{eq:approx-momentum}\\
\C^{-1}\partial_t\TT^n+\xi^n&=\eps(u_t^n)\,,\qquad \xi^n=G_n(\theta^n,\TT^n)\,,\label{eq:approx-stress}\\
\partial_t\theta^n-\Delta\theta^n+\theta^n\Div u_t^n&=q^n\,,\qquad q^n=\xi^n:\TT^n\,,\label{eq:approx-heat}
\end{align}
\end{subequations}
for all $w\in W_n$ and with $\partial_\nu\theta^n=0$ and the preceding initial data. Temperature is not spatially discretized and $q^n$ is not truncated. The only approximation parameter is $n$.

\begin{proposition}[Solvability at a fixed approximation level]\label{prop:fixed-level}
For each $n$, system~\eqref{eq:approx} has a solution on $[0,T]$ with
\begin{equation}\label{eq:fixed-regularity}
\begin{gathered}
u^n\in W^{2,2}(0,T;W_n)\,,\qquad \TT^n,\partial_t\TT^n,\xi^n,q^n\in L^\infty(Q)\,,\\
\TT^n\in W^{1,p}(0,T;L^p(\Omega;\Ssym))\quad(1\le p<\infty)\,,\\
\theta^n\in W^{1,p}(0,T;L^p(\Omega))\cap L^p(0,T;W^{2,p}(\Omega))
\quad(1<p<\infty)\,,\\
0<\underline\theta^n\le\theta^n\le\overline\theta^n<\infty\quad\text{a.e. in }Q\,.
\end{gathered}
\end{equation}
These bounds may depend on $n$. The exact energy and logarithmic entropy identities \eqref{eq:energy-formal} and \eqref{eq:free-energy-formal} hold with the approximate fields and $\ip{\xi^n}{\TT^n}=q^n$.
\end{proposition}

\begin{proof}
Fix $n$ and initially a time interval $I_h=[0,h]$. All constants in the local construction may depend on $n$ and the bounded initial data.

\noindent
\emph{The stress evolution.} Given $v\in C^1(I_h;W_n)$ and $\vartheta\in L^2(I_h\times\Omega)$ with $\vartheta\ge0$, solve for almost every $x$ the Carath\'eodory equation
\begin{equation}\label{eq:frozen-stress}
\C^{-1}(x)\partial_t S+G_n(\vartheta,S)=\eps(v_t)\,,\qquad S(0,x)=\TT_0^n(x)\,.
\end{equation}  
For almost every fixed \(x\in\Omega\), Lemma~\ref{lem:geometry} gives measurability of the right-hand side with respect to \(t\) and local Lipschitz continuity with respect to \(S\in\Ssym\), uniformly on bounded subsets of \(\Ssym\). Hence the pointwise problem admits a unique local solution. Measurable dependence on \(x\) follows, for instance, from the Picard construction applied to truncated vector fields. Testing \eqref{eq:frozen-stress} by \(S(t,x)\) yields, for almost every \((t,x)\in I_h\times\Omega\),
$$
\frac12\partial_t
\left(\ip{\C^{-1}(x)S(t,x)}{S(t,x)}\right)=\ip{\eps(v_t)(t,x)}{S(t,x)} -
\ip{G_n\bigl(\vartheta(t,x),S(t,x)\bigr)}{S(t,x)}\,.
$$
By \eqref{eq:work-identity}, the last term is nonnegative, and therefore
$$
\frac12\partial_t\left(\ip{\C^{-1}(x)S}{S}\right)\le |\eps(v_t)|\,|S| \qquad\text{a.e. in }I_h\times\Omega\,.$$
Using the uniform ellipticity of \(\C^{-1}(x)\), applying the preceding inequality to a regularized square root of \(\ip{\C^{-1}(x)S}{S}\) and then letting the regularization parameter tend to zero, we obtain
\begin{equation}\label{eq:frozen-stress-bound}
|S(t,x)|\le C_{\C}\left(|\TT_0^n(x)|+\int_0^t |\eps(v_t)(s,x)|\dd s \right)
\end{equation}
for almost every \((t,x)\in I_h\times\Omega\), where \(C_{\C}>0\) depends only on the ellipticity constants in \eqref{eq:elasticity}. The estimate \eqref{eq:frozen-stress-bound} excludes finite-time blow-up and therefore extends the pointwise stress solution to the whole interval \(I_h\).  On a bounded $C^1(I_h;W_n)$ set, $S$ is bounded in $ L^\infty(\Omega;\Ssym)$, independently of $\vartheta$. Therefore
\begin{equation}\label{eq:frozen-source-bound}
q(\vartheta,S):=G_n(\vartheta,S):S\ge0\,,\qquad
\norm{q(\vartheta,S)}_{L^\infty(I_h\times\Omega)}+\norm{\partial_tS}_{L^\infty(I_h\times\Omega)}\le C_{n,R}\,.
\end{equation}
No spatial derivative of $S$ or of $\vartheta$ is required here.

\noindent
\emph{The heat and Galerkin equations.} For this $S$, let $\Theta$ solve
\begin{equation}\label{eq:frozen-heat}
\Theta_t-\Delta\Theta+(\Div v_t)\Theta=q(\vartheta,S),\qquad
\partial_\nu\Theta=0,\qquad \Theta(0)=\theta_{0,n}.
\end{equation}
The coefficient $\Div v_t$ is bounded. Hence standard linear parabolic theory for the homogeneous Neumann problem yields a unique solution $ \Theta\in L^\infty\bigl(I_h;L^2(\Omega)\bigr) \cap
L^2\bigl(I_h;H^1(\Omega)\bigr)$, with time derivative in $\partial_t\Theta\in
L^2\bigl(I_h;(H^1(\Omega))^*\bigr)$, see, for instance, \cite{LadyzhenskayaSolonnikovUraltseva1968}. Comparison, justified by testing the positive and negative parts of differences, gives for $B=\norm{\Div v_t}_{L^\infty(I_h\times\Omega)}$ and $Q_* =\norm{q(\vartheta,S)}_{L^\infty(I_h\times\Omega)}$,
\begin{equation}\label{eq:heat-comparison}
\bigl(\operatorname*{ess\,inf}_{x\in\Omega}\theta_{0,n}(x)\bigr)e^{-Bt} \le \Theta(t,x) \le
e^{Bt} \left( \norm{\theta_{0,n}}_{L^\infty(\Omega)} +tQ_* \right)
\end{equation}
for almost every \((t,x)\in I_h\times\Omega\). In particular, $\Theta$ is positive. The bounded right-hand side $q-(\Div v_t)\Theta$ and the smooth compatible initial data give standard \(L^p\)-regularity for the heat equation \cite{LadyzhenskayaSolonnikovUraltseva1968} further gives $ \Theta\in W^{1,p}\bigl(I_h;L^p(\Omega)\bigr) \cap L^p\bigl(I_h;W^{2,p}(\Omega)\bigr)$ for every $1<p<\infty$.
 Finally define $U\in W^{2,2}(I_h;W_n)$ as the solution of
\begin{equation}\label{eq:frozen-Galerkin}
\int_\Omega \ip{U_{tt}}{w}\dd x=\int_\Omega \ip{f}{w}\dd x-\int_\Omega \ip{S}{\eps(w)}\dd x +
\int_\Omega \Theta\Div w\dd x
\end{equation}
for every \(w\in W_n\) and for almost every \(t\in I_h\), together with $U(0)=u_0^n$ and $ U_t(0)=u_1^n$. This is a forced finite-dimensional linear ODE in the Galerkin space \(W_n\). Calculations of the same general type for the thermal part of coupled inelastic systems may also be found in \cite{CieslakOwczarekWielgos2026}.

\noindent
\emph{Continuity of the coupled map.} Consider $\mathcal F(v,\vartheta)=(U,\Theta)$ on the Banach space $X_h=C^1(I_h;W_n)\times L^2(I_h\times\Omega)$ equipped with the norm 
\begin{equation}\label{eq:Xh-norm}
\|(v,\vartheta)\|_{X_h} := \|v\|_{C^1(I_h;W_n)} + \|\vartheta\|_{L^2(I_h\times\Omega)}\,,
\end{equation}
Since \(W_n\) is finite-dimensional, this choice is equivalent to the corresponding \(C^1\)-norm induced by any other norm on \(W_n\). In particular, bounded subsets of \(C^1(I_h;W_n)\) are bounded in the spatial norms required below, with constants depending on \(n\). Suppose $(v_j,\vartheta_j)\to(v,\vartheta)$ in $X_h$ inside a bounded set. The corresponding stresses share the bound \eqref{eq:frozen-stress-bound}. For \(D_j=S_j-S\), subtraction of the corresponding stress equations, monotonicity of \(G_n(\vartheta_j,\cdot)\), and Young's inequality yield
\begin{equation}
\label{eq:stress-continuity}
\begin{aligned}
\frac{\mathrm d}{\mathrm dt} \int_\Omega \ip{\C^{-1}(x)D_j}{D_j}\dd x &\le C
\norm{D_j}_{L^2(\Omega)}^2 + C\norm{\eps(v_{j,t}-v_t)}_{L^2(\Omega)}^2\\
&\quad +C \norm{G_n\bigl(\vartheta_j,S\bigr) - G_n\bigl(\vartheta,S\bigr) }_{L^2(\Omega)}^2\,.
\end{aligned}
\end{equation}
The last term converges to zero in \(L^1(I_h)\). Indeed, $\vartheta_j\to\vartheta$ in measure on $I_h\times\Omega$, the map $(\theta,S)\longmapsto G_n(\theta,S)$ is continuous, and the common bound for the stress fields \(S_j\) and \(S\) provides an \(L^\infty(I_h\times\Omega;\Sdev)\)-dominating bound for $G_n(\vartheta_j,S)-G_n(\vartheta,S)$. Gronwall's inequality therefore gives $S_j\to S$ strongly in $C\bigl(I_h;L^2(\Omega;\Ssym)\bigr)$. The same boundedness and continuity properties imply $q(\vartheta_j,S_j) \to q(\vartheta,S)$ strongly in $L^2(I_h\times\Omega)$. Subtracting the corresponding heat equations and using the standard \(L^2(\Omega)\)-energy estimate, together with $\Div v_{j,t}\to\Div v_t$ uniformly on $I_h\times\Omega$ yields
$$
\Theta_j\to\Theta \qquad\text{strongly in } C\bigl(I_h;L^2(\Omega)\bigr) \cap L^2\bigl(I_h;H^1(\Omega)\bigr)\,.
$$
Finally, equation~\eqref{eq:frozen-Galerkin} gives $U_j\to U$ strongly in $C^1(I_h;W_n)$. Continuity of \(r\) therefore suffices for continuity of the coupled map \(\mathcal F:X_h\to X_h\).

\noindent
\emph{Compactness and invariance.} Fix \(R>0\) and consider inputs \((v,\vartheta)\in X_h\) satisfying
$$
\norm{v}_{C^1(I_h;W_n)}+\norm{\vartheta}_{L^2(I_h\times\Omega)}\le R\,, \qquad \vartheta\ge0 \quad\text{a.e. in }I_h\times\Omega\,.
$$
The estimates \eqref{eq:frozen-stress-bound} and \eqref{eq:frozen-source-bound}, together with the finite dimensionality of \(W_n\), provide bounds for \(S\) and \(q(\vartheta,S)\) depending only on
\(n\), \(R\), and the fixed initial data. In particular, the heat energy estimates bound $\Theta$ in $L^2(I_h;H^1(\Omega))$ and $\partial_t\Theta$ in $L^2(I_h;(H^1(\Omega))^*)$. The compact embedding $H^1(\Omega)\Subset L^2(\Omega)$ and time compactness \cite{Simon} make the temperature outputs precompact in $L^2(I_h\times\Omega)$. On the other hand, equation~\eqref{eq:frozen-Galerkin} yields bounds $U$ in $W^{2,2}(I_h;W_n)$, compactly embedded into $C^1(I_h;W_n)$ since $W_n$ is finite-dimensional. Consequently, the map $\mathcal F:X_h\to X_h$ such that $\mathcal F(v,\vartheta)=(U,\Theta)$, maps bounded subsets of \(X_h\) into relatively compact subsets of \(X_h\).

\noindent
We next construct a closed convex set that is invariant under \(\mathcal F\).
Choose constants \(R_v>0\) and \(R_\vartheta>0\) such that $R_v> \norm{u_0^n+t u_1^n}_{C^1(I_h;W_n)}$ and set
$$
\mathcal K_h
=
\left\{
(v,\vartheta)\in X_h:
\begin{array}{l}
v(0)=u_0^n\,,\quad v_t(0)=u_1^n\,,\\[1mm]
\norm{v}_{C^1(I_h;W_n)}\le R_v\,,\\[1mm]
\vartheta\ge0\ \text{a.e. in }I_h\times\Omega\,,\\[1mm]
\norm{\vartheta}_{L^2(I_h\times\Omega)}\le R_\vartheta
\end{array}
\right\}\,.
$$
Then \(\mathcal K_h\) is a nonempty closed convex subset of \(X_h\). For \((v,\vartheta)\in\mathcal K_h\), the comparison estimate \eqref{eq:heat-comparison} and the bounds obtained above imply $\norm{\Theta}_{L^\infty(I_h\times\Omega)} \le C_{n,R_v,R_\vartheta}$ and therefore
$$
\norm{\Theta}_{L^2(I_h\times\Omega)} \le C_{n,R_v,R_\vartheta}|I_h\times\Omega|^{1/2} \le C_{n,R_v,R_\vartheta}\sqrt h\,.
$$
Moreover, \eqref{eq:heat-comparison} gives $\Theta\ge0$ a.e. in $I_h\times\Omega$. To control the displacement component, integrate \eqref{eq:frozen-Galerkin} in time. Since \(W_n\) is finite-dimensional, all spatial norms on \(W_n\) are equivalent, and the stress and temperature terms are bounded by the estimates above. The forcing term satisfies
$$
\int_0^h \norm{f(t)}_{L^2(\Omega;\R^3)}\dd t\le\sqrt h\, \norm{f}_{L^2(I_h\times\Omega;\R^3)}\,.
$$
It follows that $\norm{U_t-u_1^n}_{C(I_h;W_n)} \le C_{n,R_v,R_\vartheta}\sqrt h$ and after one further integration in time
$$
\norm{U-u_0^n-tu_1^n}_{C(I_h;W_n)} \le C_{n,R_v,R_\vartheta}h^{3/2}\,.
$$
The constants in these estimates remain bounded as \(h\downarrow0\), while the right-hand sides tend to zero. Consequently, after first fixing \(R_v\) and \(R_\vartheta\), one may choose
\(h>0\) sufficiently small so that $\mathcal F(\mathcal K_h)\subset\mathcal K_h$. The map \(\mathcal F\) is continuous by the preceding argument and has relatively compact range in \(X_h\). Hence Schauder's fixed-point theorem provides $(u^n,\theta^n)\in\mathcal K_h$ such that $\mathcal F(u^n,\theta^n)=(u^n,\theta^n)$. Together with the associated stress \(\TT^n\), this yields a solution of \eqref{eq:approx} on \(I_h\). Only existence is needed here; uniqueness of the coupled fixed-point problem is not required.

\noindent
\emph{Energy and continuation.} At a fixed point the regularity already obtained justifies testing \eqref{eq:approx-momentum} by $u_t^n$ and \eqref{eq:approx-stress} by $\TT^n$. Adding the integrated heat equation gives the exact energy identity
\begin{equation}\label{eq:approx-energy}
\begin{aligned}
E_n(t)&=E_n(0)+\int_0^t\int_\Omega\ip{f}{u_t^n}\dd x\dd s\,\\
E_n(t) &= \frac12 \norm{u_t^n(t)}_{L^2(\Omega)}^2 + \frac12 \int_\Omega \ip{\C^{-1}\TT^n}{\TT^n}\dd x + \int_\Omega \theta^n\dd x \,.
\end{aligned}
\end{equation}
Young's inequality and Gronwall give a bound for $E_n$ on every local solution up to time $T$, independent of its interval of existence. Finite dimensionality therefore bounds $u_t^n$ in $L^\infty(0,T;W^{1,\infty}(\Omega;\R^3))$, with an $n$-dependent constant. Equation~\eqref{eq:frozen-stress-bound} bounds $\TT^n$ pointwise on the same interval. Hence \eqref{eq:frozen-source-bound} bounds $q^n$ and $\partial_t\TT^n$, and \eqref{eq:heat-comparison} gives uniform-in-time upper and positive lower bounds for $\theta^n$ at this fixed $n$ i.e. $\underline\theta_n \le \theta^n(t,x) \le \overline\theta_n$ for a.e. $(t,x)\in Q$. At this stage the constants may depend on \(n\), but they are independent of the endpoint of the particular local solution under consideration.

\noindent
These bounds give a common local existence length for all restart data encountered before $T$: the fixed-point construction at a restart time \(t_0\in[0,T)\) requires only bounds for
$$
u^n(t_0),\,u_t^n(t_0)\in W_n\,, \qquad \TT^n(t_0)\in L^\infty(\Omega;\Ssym)\,, \qquad \theta^n(t_0)\in L^\infty(\Omega)
$$
together with a positive lower bound for \(\theta^n(t_0)\) and the fixed force norm $\norm{f}_{L^2(Q;\mathbb R^3)}$. Since \(W_n\) is finite-dimensional, bounds in \(W_n\) also control the spatial norms of \(u^n(t_0)\) and \(u_t^n(t_0)\) required in the local construction, with constants depending on \(n\). The linear heat problem can be restarted from any datum $\theta^n(t_0)\in L^\infty(\Omega)$ such that $\theta^n(t_0)\ge \underline\theta_n>0$ since the local energy and compactness estimates require no spatial derivatives of the restart temperature. Together with the bounds for $u^n(t_0),\,u_t^n(t_0)\in W_n$ and $\TT^n(t_0)\in L^\infty(\Omega;\Ssym)$ this yields a restart interval of length bounded from below independently of \(t_0\). The pieces match at the restart times, and after finitely many restarts the solution is defined on \([0,T]\).

\noindent
On the concatenated interval the heat equation has a bounded right-hand side and the original smooth initial datum, so global finite-$p$ parabolic regularity follows from the linear equation. The stress equation supplies $\partial_t\TT^n\in L^\infty(Q)$ and thus the asserted $W^{1,p}\bigl(0,T;L^p(\Omega;\Ssym)\bigr)$ regularity for finite $p$. Positivity permits division of the heat equation by $\theta^n$ and gives
\begin{equation}\label{eq:approx-log}
\partial_t(\tau^n+\Div u^n)-\Delta\tau^n
=\frac{q^n}{\theta^n}+|\nabla\tau^n|^2,\qquad\tau^n=\log\theta^n.
\end{equation}
The chain rule is legitimate at this fixed level because $\theta^n$ and its reciprocal are bounded and the parabolic regularity holds for every finite $p$. Integration of \eqref{eq:approx-log}, followed by subtraction from \eqref{eq:approx-energy}, proves the remaining identity.
\end{proof}

\noindent
A related semi-Galerkin construction was used in \cite{CieslakOwczarekWielgos2026}, where separate approximation levels were introduced for the displacement and the stress tensor. We also refer there for more details on the standard fixed-point and semi-Galerkin arguments. Here only the displacement is Galerkin-discretized, so that the limiting procedure involves a single approximation parameter, as in the three-dimensional thermoelastic construction of  \cite{CieslakMuhaTrifunovic2024}.
\section{Uniform estimates}
\label{sec:estimates}

All constants in this section are independent of $n$, unless a dependence is explicitly indicated. The energy convergence of the initial data follows from \eqref{eq:initial-temp-approx}, \eqref{eq:initial-stress-approx}, and the Galerkin projections.
\begin{lemma}[Energy, entropy, and volumetric control]\label{lem:uniform}
The approximate solutions satisfy
\begin{equation}\label{eq:uniform-basic}
\begin{aligned}
&\norm{u_t^n}_{L^\infty(0,T;L^2(\Omega))} +\norm{\TT^n}_{L^\infty(0,T;L^2(\Omega))} +\norm{\theta^n}_{L^\infty(0,T;L^1(\Omega))}\le C\,,\\
&\norm{u^n}_{L^\infty(0,T;L^2(\Omega))} +\norm{\Div u^n}_{L^\infty(0,T;L^2(\Omega))}\le C\,,\\
&\norm{\tau^n}_{L^\infty(0,T;L^1(\Omega))} +\norm{\tau^n}_{L^2(0,T;H^1(\Omega))}
+\norm{\nabla\tau^n}_{L^2(Q)} +\int_Q\frac{q^n}{\theta^n}\dd x\dd t\le C\,.
\end{aligned}
\end{equation}
For every finite $s>3$,
\begin{equation}\label{eq:tau-time}
\norm{\partial_t\tau^n}_{L^1(0,T;Y_s)}\le C_s\,.
\end{equation}
\end{lemma}

\begin{proof}
From \eqref{eq:approx-energy}, ellipticity, and $\int_\Omega\ip{f}{u_t^n}\dd x \le\frac12\norm{f}_{L^2(\Omega)}^2 +\frac12\norm{u_t^n}_{L^2(\Omega)}^2$, Gronwall gives the first line of \eqref{eq:uniform-basic}. Integrating velocity controls $u^n$. Taking the trace of \eqref{eq:approx-stress} and using $\tr\xi^n=0$ gives the exact identity
\begin{equation}\label{eq:approx-volume}
\Div u^n(t)=\Div u_0^n
+\tr\bigl(\C^{-1}(\TT^n(t)-\TT_0^n)\bigr)\,.
\end{equation}
This proves the volumetric bound, including for heterogeneous anisotropic $\C$.

\noindent
Integration of \eqref{eq:approx-log} uses the Neumann condition and
$\int_\Omega\Div u^n\dd x=0$, and yields
\begin{equation}\label{eq:integrated-entropy}
\int_\Omega\tau^n(t)\dd x-\int_\Omega\log\theta_{0,n}\dd x = \int_0^t\!\int_\Omega \left(\frac{q^n}{\theta^n}+|\nabla\tau^n|^2\right)\dd x\dd s\,.
\end{equation}
The right-hand side is nonnegative. Also $\{\log s\}_+\le s$ for $s>0$. The energy bound and the lower bound for $\int_\Omega\tau^n\dd x$ therefore bound both positive and negative parts of $\tau^n$ in $L^\infty(0,T;L^1(\Omega))$. They also bound the right-hand side of \eqref{eq:integrated-entropy}. Poincar\'e's inequality with the spatial mean gives
\[
\norm{\tau^n(t)}_{H^1(\Omega)}
\le C\left( \norm{\nabla\tau^n(t)}_{L^2(\Omega)} +\norm{\tau^n(t)}_{L^1(\Omega)} \right),
\]
which proves the remaining spatial estimate. Equivalently, subtraction of \eqref{eq:integrated-entropy} from \eqref{eq:approx-energy} gives
\begin{equation}\label{eq:approx-combined}
\begin{aligned}
E_n(t)-\int_\Omega\tau^n(t)\dd x &+\int_0^t\!\int_\Omega \left(\frac{q^n}{\theta^n}+|\nabla\tau^n|^2\right)\dd x\dd s\\
&=E_n(0)-\int_\Omega\log\theta_{0,n}\dd x +\int_0^t\!\int_\Omega\ip{f}{u_t^n}\dd x\dd s\,.
\end{aligned}
\end{equation}
To estimate time derivatives, test \eqref{eq:approx-log} by $\zeta\in W^{1,s}(\Omega)$. Since $s>3$, $W^{1,s}(\Omega)\hookrightarrow L^\infty(\Omega)$ and $W^{1,s}(\Omega)\hookrightarrow H^1(\Omega)$. The Neumann Laplacian is used variationally, through $-\int_\Omega\ip{\nabla\tau^n}{\nabla\zeta}\dd x$, Furthermore, $u_t^n=0$ on the boundary, so $\dual{\Div u_t^n}{\zeta}
=-\int_\Omega\ip{u_t^n}{\nabla\zeta}\dd x$. Consequently,
\[
\norm{\partial_t\tau^n(t)}_{Y_s} \le C_s\left( \norm{\nabla\tau^n(t)}_{L^2(\Omega)} +\norm{u_t^n(t)}_{L^2(\Omega)} +\norm{q^n(t)/\theta^n(t)}_{L^1(\Omega)} +\norm{\nabla\tau^n(t)}_{L^2(\Omega)}^2 \right)\,.
\]
Integration over $(0,T)$ proves \eqref{eq:tau-time}. We use $Y_s=(W^{1,s}(\Omega))^*$ rather than $(W^{1,s}_0(\Omega))^*$, since the Neumann boundary condition allows test functions in $W^{1,s}(\Omega)$ without any restriction on their boundary trace.
\end{proof}

\noindent
The preceding energy--entropy estimates do not provide an unweighted $L^1$ bound for the plastic work $q^n$. The next two results show that the available entropy control is nevertheless sufficient to obtain a quantitative estimate of the overstress, that is, of the distance of $\dev\TT^n$ from the moving yield set $r(\theta^n)C_0$. This estimate is one of the key compactness mechanisms of the approximation: it yields the limiting stress admissibility and will later enter the construction of the weighted constitutive measures and their lifts.

\begin{lemma}[Square-root plastic work]\label{lem:square-root}
One has
\begin{equation}\label{eq:square-root}
\int_0^T\!\int_\Omega \sqrt{q^n}\dd x\dd t\le C\,,\qquad \int_0^T\!\int_\Omega d_n(1+d_n^2)^{(n-1)/4}\dd x\dd t\le C\,.
\end{equation}
\end{lemma}

\begin{proof}
Positivity of the temperature gives
\[
\int_0^T\!\int_\Omega\sqrt{q^n}\dd x\dd t = \int_0^T\!\int_\Omega \sqrt{\frac{q^n}{\theta^n}}\sqrt{\theta^n}\,\dd x\dd t \le \left(\int_0^T\!\int_\Omega\frac{q^n}{\theta^n}\dd x\dd t\right)^{1/2} \left(\int_0^T\!\int_\Omega\theta^n\dd x\dd t\right)^{1/2} \le C\,.
\]
The second assertion follows by taking square roots in the lower bound of \eqref{eq:work-identity}. No unweighted $L^1((0,T)\times\Omega)$ estimate for $q^n$ is used.
\end{proof}

\begin{proposition}[Quantitative convergence to the moving constraint]\label{prop:overstress}
For every $p>1$ and $n>2p-1$,
\begin{equation}\label{eq:algebraic-rate}
\int_0^T\!\int_\Omega d_n^p\dd x\dd t\le Cc_{n,p}\,,\qquad c_{n,p} = \left(\frac{2(p-1)}{n-2p+1}\right)^{(p-1)/2} \left(\frac{n-2p+1}{n-1}\right)^{(n-1)/4}\,.
\end{equation}
In particular, $c_{n,p}\le C_p n^{-(p-1)/2}$ for all sufficiently large $n$. For each fixed $1\le p<\infty$, the sharper bound \eqref{eq:main-rate} holds. Moreover, for every $\eta>0$,
\begin{equation}\label{eq:overstress-tail}
\bigl|\{(t,x)\in Q:\ d_n(t,x)\ge\eta\}\bigr| \le\frac{C}{\eta(1+\eta^2)^{(n-1)/4}}\,.
\end{equation}
Finally,
\begin{equation}\label{eq:collapse-approx}
\norm{\widehat\TT^n-\TT^n}_{L^p(0,T;L^p(\Omega;\Ssym))} = \norm{d_n}_{L^p(0,T;L^p(\Omega))}\,,\qquad \int_{\{r_n=0\}}|\dev\TT^n|^2\dd x\dd t \le C\frac{\log n}{n}
\end{equation}
for all sufficiently large $n$ in the second inequality.
\end{proposition}

\begin{proof}
Write
\[
s^p = s(1+s^2)^{(n-1)/4} \frac{s^{p-1}}{(1+s^2)^{(n-1)/4}} \qquad (s\ge0)\,.
\]
For $p>1$, set
\[
g_{n,p}(s)=\frac{s^{p-1}}{(1+s^2)^{(n-1)/4}}\,.
\]
If $n>2p-1$, then
\[
\frac{\dd}{\dd s}\log g_{n,p}(s) = \frac{p-1}{s}-\frac{n-1}{2}\frac{s}{1+s^2}\,,
\]
and therefore the unique positive critical point is determined by $s^2=\frac{2(p-1)}{n-2p+1}$. The function $g_{n,p}$ increases before this point and decreases afterwards, so this critical point gives its maximum. Substitution yields
\[
\sup_{s\ge0}g_{n,p}(s) = \left(\frac{2(p-1)}{n-2p+1}\right)^{(p-1)/2} \left(\frac{n-2p+1}{n-1}\right)^{(n-1)/4} = c_{n,p}\,.
\]
Consequently, $s^p\le c_{n,p}\,s(1+s^2)^{(n-1)/4}$ for $s\ge0$. Taking $s=d_n(t,x)$, integrating over $(0,T)\times\Omega$, and using Lemma~\ref{lem:square-root}, we obtain
\[
\int_0^T\!\int_\Omega d_n^p\dd x\dd t \le c_{n,p} \int_0^T\!\int_\Omega d_n(1+d_n^2)^{(n-1)/4}\dd x\dd t \le Cc_{n,p}\,,
\]
which proves \eqref{eq:algebraic-rate}. For sufficiently large $n$, the first factor in the definition of $c_{n,p}$ is bounded by $C_p n^{-(p-1)/2}$, while the second factor is at most one. Hence $c_{n,p}\le C_p n^{-(p-1)/2}$. The map $s\longmapsto s(1+s^2)^{(n-1)/4}$ is increasing on $[0,\infty)$. Thus, on the set $\{d_n\ge\eta\}$
\[
d_n(1+d_n^2)^{(n-1)/4} \ge \eta(1+\eta^2)^{(n-1)/4}\,.
\]
Using again Lemma~\ref{lem:square-root}, we therefore obtain
\[
C \ge \int_{\{d_n\ge\eta\}} d_n(1+d_n^2)^{(n-1)/4}\dd x\dd t \ge \eta(1+\eta^2)^{(n-1)/4} |\{d_n\ge\eta\}|\,,
\]
which proves \eqref{eq:overstress-tail}. For the logarithmic improvement, fix $p\ge1$ and let
\[
\beta_n=\left(\frac{16\log n}{n-1}\right)^{1/2}\,.
\]
For all sufficiently large $n$, $\beta_n^2\le1$, and hence $\log(1+\beta_n^2)\ge\frac{\beta_n^2}{2}$. It follows that
\[
\frac{n-1}{4}\log(1+\beta_n^2) \ge \frac{n-1}{8}\beta_n^2 = 2\log n\,,
\]
and therefore $(1+\beta_n^2)^{(n-1)/4}\ge n^2$. For $p=1$, the function $s\longmapsto (1+s^2)^{-(n-1)/4}$ is decreasing. For $p>1$, the preceding computation shows that $s\longmapsto s^{p-1}(1+s^2)^{-(n-1)/4}$ is decreasing once $s$ lies to the right of its critical point. Since, for fixed $p$ and sufficiently large $n$
\[
\beta_n^2 = \frac{16\log n}{n-1} > \frac{2(p-1)}{n-2p+1}\,,
\]
the function is decreasing on $[\beta_n,\infty)$. We now split the integral into the regions $\{d_n<\beta_n\}$ and $\{d_n\ge\beta_n\}$. On the first one,
\[
\int_{\{d_n<\beta_n\}}d_n^p\dd x\dd t \le |Q|\beta_n^p\,.
\]
On the second one, using the monotonicity just established
\[
\frac{d_n^{p-1}}{(1+d_n^2)^{(n-1)/4}} \le \frac{\beta_n^{p-1}}{(1+\beta_n^2)^{(n-1)/4}}\,.
\]
Hence, by Lemma~\ref{lem:square-root},
\[
\begin{aligned}
\int_{\{d_n\ge\beta_n\}}d_n^p\dd x\dd t &\le \frac{\beta_n^{p-1}}{(1+\beta_n^2)^{(n-1)/4}}
\int_0^T\!\int_\Omega d_n(1+d_n^2)^{(n-1)/4}\dd x\dd t\\
&\le C\frac{\beta_n^{p-1}}{(1+\beta_n^2)^{(n-1)/4}} \le C\beta_n^{p-1}n^{-2}\,.
\end{aligned}
\]
Combining the two regions gives, then taking the $p$th root we get\eqref{eq:main-rate}. Finally, the projection changes only the deviatoric part, and therefore $|\widehat\TT^n-\TT^n|=d_n$. On $\{r_n=0\}$ one has $r_nC_0=\{0\}$ and hence $Z_n=0$. Consequently, $|\dev\TT^n|=d_n$ on $\{r_n=0\}$. Applying \eqref{eq:main-rate} with $p=2$ gives
\[
\int_{\{r_n=0\}}|\dev\TT^n|^2\dd x\dd t
\le
\int_0^T\!\int_\Omega d_n^2\dd x\dd t
\le
C\frac{\log n}{n},
\]
which proves \eqref{eq:collapse-approx}.
\end{proof}

\noindent
The preceding estimates are a key quantitative consequence of the Norton--Hoff penalization. Although the available entropy bound controls only the square root of the plastic work, the rapidly increasing factor $(1+d_n^2)^{(n-1)/4}$ forces the overstress $d_n$ to zero strongly in every finite $L^p(Q)$. In particular, the choice $\beta_n\simeq\sqrt{\log n/n}$ extracts from this penalization the rate required below to identify the limiting moving yield constraint, while \eqref{eq:overstress-tail} shows that regions with a fixed positive overstress become rapidly negligible.

\section{Compactness and passage to the limit}
\label{sec:limit}

\subsection{Logarithmic temperature and thermal energy}

\begin{proposition}[Thermal compactness]\label{prop:thermal}
After extraction, the convergences \eqref{eq:main-tau} hold, together with
\begin{equation}\label{eq:tau-weak}
\tau^n\rightharpoonup\tau
\quad\text{in }L^2(0,T;H^1(\Omega))\,,\qquad \tau\in L^\infty(0,T;L^1(\Omega))\cap BV([0,T];Y_s)\,.
\end{equation}
The thermal measures satisfy \eqref{eq:thermal-split} and $e^\tau\in L^\infty(0,T;L^1(\Omega))$.
\end{proposition}

\begin{proof}
For $1\le q<6$, Rellich's theorem gives $H^1(\Omega)\Subset L^q(\Omega)$ and $L^q(\Omega)\hookrightarrow Y_s$ continuously. The estimates \eqref{eq:uniform-basic} and \eqref{eq:tau-time} satisfy the finite-time-exponent compactness criterion of Aubin--Lions--Simon \cite{Simon}, and therefore $\tau^n\longrightarrow\tau$ strongly in $L^2(0,T;L^q(\Omega))$ for every $1\le q<6$. A diagonal extraction for exponents approaching $6$, followed by interpolation or by the continuous inclusions on $\Omega$, gives one subsequence for all $q<6$. In particular, after a further subsequence, $\tau^n\longrightarrow\tau$ a.e. in $Q$. Fatou's lemma on time slices gives $\tau\in L^\infty(0,T;L^1(\Omega))$. The derivative estimate \eqref{eq:tau-time} passes to a $Y_s$-valued finite-variation derivative. Since $Y_s$ is reflexive for finite $s$, the uniform variation bound together with convergence in $L^1(0,T;Y_s)$ yields the asserted representative $\tau\in BV([0,T];Y_s)$. A countable sequence of exponents $s\downarrow3$ suffices to obtain the conclusion for every finite $s>3$ by continuous embeddings.

\noindent
Interpolation in space gives
\begin{equation}\label{eq:tau-critical-integrability}
\int_0^T \norm{\tau^n(t)}_{L^{8/3}(\Omega)}^{8/3}\dd t \le C \norm{\tau^n}_{L^\infty(0,T;L^1(\Omega))}^{2/3} \norm{\tau^n}_{L^2(0,T;L^6(\Omega))}^{2} \le C\,.
\end{equation}
Strong convergence in $L^2(Q)$ together with this bound implies $\tau^n\longrightarrow\tau$ strongly in $L^p(Q)$ for every $1\le p<\frac83$. Exponentiation gives $\theta^n\longrightarrow e^\tau>0$ a.e. in $Q$ and Fatou's lemma again yields $e^\tau\in L^\infty(0,T;L^1(\Omega))$. Since $r$ is bounded and continuous, $r_n=r(\theta^n)\longrightarrow r(e^\tau)$ a.e. in $Q$ and by dominated convergence 
\begin{equation}\label{eq:srtong_rn}
r_n\longrightarrow r(e^\tau) \quad\text{strongly in }L^p(Q) \qquad\text{for every }1\le p<\infty\,.
\end{equation}
The uniform bound $\norm{\theta^n}_{L^\infty(0,T;L^1(\Omega))}\le C$ gives, by weak-star compactness in the dual of the separable space $L^1(0,T;C(\overline\Omega))$, a positive limit $\Theta\in L^\infty_{w*}(0,T;\M(\overline\Omega))$. For nonnegative $\zeta\in C(\overline\Omega)$ and $\psi\in L^\infty(0,T)$, Fatou's lemma gives
\[
\int_0^T\psi(t)\int_\Omega \zeta e^{\tau}\dd x\dd t \le \liminf_{n\to\infty} \int_0^T\!\int_\Omega\psi(t)\zeta\theta^n\dd x\dd t = \int_0^T\psi(t)\dual{\Theta(t)}{\zeta}\dd t\,.
\]
A countable dense set of nonnegative functions $\zeta\in C(\overline\Omega)$ and the arbitrariness of nonnegative $\psi\in L^\infty(0,T)$ imply
\[
\Theta(t)\ge e^{\tau(t)}\dd x \qquad\text{for a.e. }t\in(0,T)\,.
\]
Their difference defines the weak-star measurable nonnegative measure $\lambda_t$ appearing in \eqref{eq:thermal-split}. In spacetime, the measure $\Theta(t)\dd t$ is absolutely continuous with respect to time. The excess $\lambda_t$ need not be singular with respect to spatial Lebesgue measure: concentration may move with $n$ and therefore have a diffuse weak limit.
\end{proof}

\subsection{Mechanical compactness and admissibility}

The uniform estimates yield, along the same subsequence
\begin{equation}\label{eq:mechanical-limits}
\begin{aligned}
u^n&\weakstar u &&\text{in }W^{1,\infty}(0,T;L^2(\Omega;\R^3))\,,\\
u_t^n&\weakstar u_t &&\text{in }L^\infty(0,T;L^2(\Omega;\R^3))\,,\\
\TT^n&\weakstar\TT &&\text{in }L^\infty(0,T;L^2(\Omega;\Ssym))\,,\\
\Div u^n&\weakstar\Div u &&\text{in }L^\infty(0,T;L^2(\Omega))\,.
\end{aligned}
\end{equation}
The weak-star notation in $W^{1,\infty}(0,T;L^2(\Omega;\mathbb R^3))$ means weak-star convergence of both the function and its time derivative in the corresponding $L^\infty$ spaces. Weak closedness of the zero-normal-trace space yields $u\in L^\infty(0,T;\Hzdiv)$. Passing to the identity $u^n(t)=u_0^n+\int_0^t u_t^n(s)\dd s$ and testing against arbitrary elements of $L^2(\Omega;\mathbb R^3)$ gives $u(t)=u_0+\int_0^t u_t(s)\,\dd s$ in $L^2(\Omega;\mathbb R^3)$. Hence $u\in C([0,T];L^2(\Omega;\mathbb R^3))$ and in particular $u(0)=u_0$ strongly in $L^2(\Omega;\mathbb R^3)$.
\begin{proposition}[Admissibility and deviatoric compactness]\label{prop:admissibility}
The assertions \eqref{eq:main-admissibility} hold. In addition,
\begin{equation}\label{eq:projected-convergence}
Z_n\weakstar\dev\TT\quad\text{in }L^\infty(Q;\Sdev),\qquad
\dev\TT^n\rightharpoonup\dev\TT\quad\text{in }L^p(Q;\Sdev),\quad 1\le p<\infty\,.
\end{equation}
\end{proposition}

\begin{proof}
Since $Z_n\in r_nC_0$ and $C_0\subset\overline B_{\kappa_0}(0)$, we have $|Z_n|\le\kappa_0r_n\le\kappa_0r_{\max}$ a.e. in $(0,T)\times\Omega$. Hence, after extraction, $Z_n\weakstar Z$ in $L^\infty(0,T;L^\infty(\Omega;\Sdev))$ for some $Z$. On the other hand, Proposition~\ref{prop:overstress} gives
\[
\norm{\dev\TT^n-Z_n}_{L^2(0,T;L^2(\Omega;\Sdev))} = \norm{d_n}_{L^2(0,T;L^2(\Omega))} \longrightarrow0\, ,
\]
while $\dev\TT^n\rightharpoonup\dev\TT$ in $L^2(0,T;L^2(\Omega;\Sdev))$. It follows that $Z=\dev\TT$, and therefore
\[
Z_n\weakstar\dev\TT \quad\text{in }L^\infty(0,T;L^\infty(\Omega;\Sdev))\,.
\]
We next identify the limiting constraint. For every fixed $W\in\Sdev$, the inclusion $Z_n\in r_nC_0$ and the definition of the support function give
\[
\ip{Z_n}{W}\le r_n h_0(W)
\qquad\text{a.e. in }(0,T)\times\Omega\,.
\]
Let $\varphi\in L^\infty((0,T)\times\Omega)$ be nonnegative. Multiplying by $\varphi$ and integrating yields
\[
\int_0^T\!\int_\Omega \varphi\ip{Z_n}{W}\,\dd x\dd t \le h_0(W) \int_0^T\!\int_\Omega \varphi r_n\dd x\dd t\,.
\]
Using the weak-star convergence of $Z_n$ and the strong convergence \eqref{eq:srtong_rn} we obtain
\[
\int_0^T\!\int_\Omega
\varphi\,\ip{\dev\TT}{W}\dd x\dd t \le h_0(W) \int_0^T\!\int_\Omega \varphi\, r(e^\tau)\dd x\dd t\,.
\]
Since $\varphi$ is arbitrary and nonnegative, $\ip{\dev\TT}{W} \le r(e^\tau)h_0(W)$ a.e. in $(0,T)\times\Omega$. Taking \(W\) first in a countable dense subset \(\{W_k\}_{k\in\mathbb N}\subset\Sdev\), the preceding inequality holds outside a null set \(N_k\) depending on \(W_k\). Since the family is countable, the union \(N=\bigcup_{k=1}^\infty N_k\) still has measure zero, and hence $\ip{\dev\TT}{W_k}\le r(e^\tau)h_0(W_k)$ holds simultaneously for all \(k\) outside \(N\). By density of \(\{W_k\}\) in \(\Sdev\) and continuity of \(h_0\), we may pass to the limit and obtain
$$
\ip{\dev\TT}{W}\le r(e^\tau)h_0(W) \qquad\text{for every }W\in\Sdev
$$
outside the same null set. The support-function characterization of the closed convex set \(r(e^\tau)C_0\) then yields
$$
\dev\TT\in r(e^\tau)C_0 \qquad\text{a.e. in }(0,T)\times\Omega\,.
$$
In particular,
\[
|\dev\TT|\le\kappa_0r(e^\tau)\le\kappa_0r_{\max} \qquad\text{a.e. in }(0,T)\times\Omega\,,
\]
which gives the asserted uniform bound. Finally, $\dev\TT^n=Z_n+(\dev\TT^n-Z_n)$ and Proposition~\ref{prop:overstress} together with the weak-star convergence of $Z_n$ yield
\[
\dev\TT^n\rightharpoonup\dev\TT \quad\text{in }L^p(0,T;L^p(\Omega;\Sdev)) \qquad\text{for every }1\le p<\infty\,.
\]
It remains to consider the zero-radius region $B_0=\{(t,x)\in(0,T)\times\Omega:\ r(e^{\tau(t,x)})=0\}$. Since $|Z_n|\le\kappa_0r_n$,
\[
\mathbf1_{B_0}|\dev\TT^n| \le \mathbf1_{B_0}|Z_n|+\mathbf1_{B_0}|\dev\TT^n-Z_n|\le \kappa_0\mathbf1_{B_0}r_n+d_n\,.
\]
The strong convergence of $r_n$ to $r(e^\tau)$ implies $\mathbf1_{B_0}r_n\longrightarrow0$ strongly in $L^p(0,T;L^p(\Omega))$ for every finite $p$, while Proposition~\ref{prop:overstress} gives the same convergence for $d_n$. Consequently,
\[
\mathbf1_{B_0}\dev\TT^n\longrightarrow0 \quad\text{strongly in }L^p(0,T;L^p(\Omega;\Sdev)) \qquad\text{for every }1\le p<\infty\,.
\]
This has a simple geometric interpretation. On the limiting zero-radius region $B_0$, the admissible deviatoric set collapses to $r(e^\tau)C_0=\{0\}$. Thus the limiting deviatoric stress must vanish there, while the preceding argument gives the stronger approximation-level conclusion that $\dev\TT^n$ itself vanishes strongly on $B_0$ in every finite $L^p$ space. This proves the remaining assertion of \eqref{eq:main-admissibility}. Notice that no quantitative rate is claimed here, since the convergence of $r_n$ to $r(e^\tau)$ is not accompanied by a quantitative rate.
\end{proof}
\noindent
The compactness established above is sufficient to pass to the momentum and entropy equations and, together with the quantitative decay of the overstress \(d_n\), to identify the limiting moving yield constraint. The full plastic flow requires a separate argument: it is first identified in a fixed negative Sobolev space, while its weighted measure limits are constructed later in Section~\ref{sec:plastic}.

\subsection{Momentum and constitutive balance}
For $\varphi\in C_c^\infty([0,T)\times\Omega;\R^3)$, we use $P_n\varphi(t)$ as a test function in \eqref{eq:approx-momentum} and integrate over $(0,T)$. The mechanical and force terms converge by \eqref{eq:mechanical-limits} and by the spectral convergence of the projected test functions. For the thermal term, we write
\[
\left| \int_0^T\!\int_\Omega \theta^n\Div(P_n\varphi-\varphi)\dd x\dd t \right| \le
\norm{\theta^n}_{L^\infty(0,T;L^1(\Omega))} \norm{\Div(P_n\varphi-\varphi)}_{L^1(0,T;L^\infty(\Omega))} \longrightarrow0\,.
\]
The remaining thermal contribution converges to $\int_0^T\dual{\Theta(t)}{\Div\varphi(t)}\dd t$. Since $u_1^n\to u_1$ strongly in $L^2(\Omega;\R^3)$, we may also pass to the initial velocity term. This proves \eqref{eq:state-momentum}. For every $w\in C_c^\infty(\Omega;\R^3)$, the limiting momentum equation shows that the map $t\longmapsto\int_\Omega\ip{u_t(t)}{w}\dd x$ is absolutely continuous on $[0,T]$, since its distributional derivative is the sum of terms generated by the stress, the thermal measure, and the force, all belonging to $L^1(0,T)$. Together with the uniform bound $u_t\in L^\infty(0,T;L^2(\Omega;\R^3))$ and a density argument, this yields $u_t\in C_w([0,T];L^2(\Omega;\R^3))$ and $u_t(0)=u_1$. We next identify the full plastic flow in a fixed negative Sobolev space. Let $\Phi\in C_c^\infty((0,T)\times\Omega;\Ssym)$. Using the constitutive equation \eqref{eq:approx-stress} and integrating the stress derivative by parts in time, we obtain
\begin{equation}\label{eq:Xi-bound}
\dual{\xi^n}{\Phi} = -\int_0^T\!\int_\Omega \ip{u_t^n}{\Div\Phi}\dd x\dd t + \int_0^T\!\int_\Omega
\ip{\C^{-1}\TT^n}{\partial_t\Phi}\dd x\dd t
\end{equation}
and $|\dual{\xi^n}{\Phi}| \le C\norm{\Phi}_{H^1((0,T)\times\Omega)}$. Hence $\xi^n$ is uniformly bounded in $H^{-1}((0,T)\times\Omega;\Ssym)$. After extraction, $\xi^n\rightharpoonup\Xi$ in $H^{-1}((0,T)\times\Omega;\Ssym)$ and passage to the limit in \eqref{eq:Xi-bound} identifies
\[
\Xi = \eps(u_t)-\partial_t(\C^{-1}\TT)
\]
in the sense of distributions on $(0,T)\times\Omega$. Since every $\xi^n$ is deviatoric, the limit $\Xi$ is trace-free as well. This proves the second part of \eqref{eq:state-constitutive}.

\noindent
Allowing test functions that reach $t=0$ produces the additional initial term in \eqref{eq:initial-constitutive}. Its convergence follows from $\TT_0^n\longrightarrow\TT_0$ strongly in $L^2(\Omega;\Ssym)$. Finally, passing to the limit in \eqref{eq:approx-volume} gives \eqref{eq:state-volume}, and this identity also ensures that the initial extension of the constitutive distribution is trace-free. Thus the full plastic flow is identified in the fixed space $H^{-1}((0,T)\times\Omega;\Ssym)$.

\subsection{Entropy and the two energy inequalities}

The uniform estimates in \eqref{eq:uniform-basic} give that after extraction of a common subsequence, there exist finite nonnegative measures $\mu_{\rm pl},\mu_{\rm th}\in\M^+(\qbar)$ such that
\begin{equation}\label{eq:production-limits}
\frac{q^n}{\theta^n}\dd x\dd t\weakstar\mu_{\rm pl}\,, \qquad |\nabla\tau^n|^2\dd x\dd t\weakstar\mu_{\rm th} \quad\text{in }\M(\qbar)\,.
\end{equation}
Let $\phi\in C_c^\infty([0,T)\times\overline\Omega)$. Testing \eqref{eq:approx-log} by $\phi$ and integrating by parts in time and space, using the Neumann boundary condition, gives
\[
\begin{aligned}
-\int_0^T\!\int_\Omega(\tau^n+\Div u^n)\partial_t\phi\dd x\dd t
&+\int_0^T\!\int_\Omega\ip{\nabla\tau^n}{\nabla\phi}\dd x\dd t =
\int_0^T\!\int_\Omega \phi\left(\frac{q^n}{\theta^n}+|\nabla\tau^n|^2\right)\dd x\dd t\\
&\quad +\int_\Omega(\log\theta_{0,n}+\Div u_0^n)\phi(0)\dd x\,.
\end{aligned}
\]
The strong convergence of $\tau^n$, the weak convergence of $\Div u^n$ and $\nabla\tau^n$, the measure convergences in \eqref{eq:production-limits}, and the convergence of the initial data allow us to pass to the limit and obtain
\[
\begin{aligned}
-\int_0^T\!\int_\Omega(\tau+\Div u)\partial_t\phi\dd x\dd t
+\int_0^T\!\int_\Omega\ip{\nabla\tau}{\nabla\phi}\dd x\dd t
&= \int_{\qbar}\phi\,\dd(\mu_{\rm pl}+\mu_{\rm th})\\
&\quad +\int_\Omega(\tau_0+\Div u_0)\phi(0)\dd x\,,
\end{aligned}
\]
which is precisely \eqref{eq:state-entropy}. Moreover, for every nonnegative $\zeta\in C(\qbar)$, weak lower semicontinuity gives
\[
\int_0^T\!\int_\Omega \zeta\,|\nabla\tau|^2\dd x\dd t \le \liminf_{n\to\infty} \int_0^T\!\int_\Omega \zeta\,|\nabla\tau^n|^2\dd x\dd t = \int_{\qbar}\zeta\,\dd\mu_{\rm th}\,.
\]
Consequently, $\mu_{\rm th}\ge|\nabla\tau|^2\dd x\dd t$ as measures on $\qbar$. In particular, the entropy inequality follows from the preceding equality and the nonnegativity of the production measures. Notice that the product $\theta^n\Div u_t^n$ is never passed separately to the limit; it has already been incorporated into the logarithmic entropy equation \eqref{eq:approx-log}.

\noindent
We next pass to the energy balance. The argument is carried out after integration in time, since the convergences in \eqref{eq:mechanical-limits} do not provide weak convergence of $\TT^n(t)$ for every fixed $t$. Let $\psi\in C_c^\infty(0,T)$ be nonnegative and set $F_\psi(s)=\int_s^T\psi(t)\dd t$. Multiplying \eqref{eq:approx-energy} by $\psi(t)$ and integrating over $(0,T)$, we obtain
\[
\int_0^T\psi(t)E_n(t)\dd t = E_n(0)\int_0^T\psi(t)\dd t + \int_0^T\psi(t)\int_0^t\!\int_\Omega\ip{f}{u_t^n}\dd x\dd s\dd t\,.
\]
By Fubini's theorem,
\[
\int_0^T\psi(t)\int_0^t\!\int_\Omega\ip{f}{u_t^n}\dd x\dd s\dd t = \int_0^T\!\int_\Omega F_\psi(t)\ip{f}{u_t^n}\dd x\dd t\,.
\]
The right-hand side converges by the weak-star convergence of $u_t^n\weakstar u_t$ in $L^\infty(0,T;L^2(\Omega;\R^3))$. For the left-hand side, weighted weak lower semicontinuity treats the kinetic and elastic terms, while $\theta^n\dd x\weakstar\Theta$ in $L^\infty_{w*}(0,T;\M(\overline\Omega))$ gives convergence of the thermal contribution. Since the initial energies converge to $E_0$, we conclude that
\[
\int_0^T\psi(t)\mathcal E(t)\dd t
\le
\int_0^T\psi(t)
\left(
E_0+\int_0^t\!\int_\Omega\ip{f}{u_t}\dd x\dd s
\right)\dd t.
\]
Since this holds for every nonnegative $\psi\in C_c^\infty(0,T)$, it follows that
\eqref{eq:state-energy-ineq} holds for almost every $t\in(0,T)$.

\noindent
The same argument applies to the combined identity \eqref{eq:approx-combined}. The term containing $\tau^n$ converges strongly. For either $m_n=\frac{q^n}{\theta^n}$ or $m_n=|\nabla\tau^n|^2$ with corresponding limit measure $\mu=\mu_{\rm pl}$ or $\mu=\mu_{\rm th}$, Fubini's theorem gives
\[
\begin{aligned}
\int_0^T\psi(t)\int_0^t\!\int_\Omega m_n\dd x\dd s\dd t =\int_0^T\!\int_\Omega F_\psi(s)m_n\dd x\dd s \longrightarrow \int_{\qbar}F_\psi(s)\dd\mu(s,x)\,.
\end{aligned}
\]
Since $F_\psi$ is continuous on $[0,T]$, the measure convergence in \eqref{eq:production-limits} applies directly. Moreover, $F_\psi(0)=\int_0^T\psi(t)\dd t$, so that possible production atoms at the initial time are retained in the limit. A second application of Fubini's theorem yields
\[
\int_{\qbar}F_\psi(s)\dd\mu(s,x) = \int_0^T\psi(t)\, \mu([0,t]\times\overline\Omega)\dd t
\]
up to the choice of the upper time endpoint at the at most countably many atoms of $\mu$. This distinction is irrelevant for an identity integrated in $t$. Passing to the limit in \eqref{eq:approx-combined} therefore gives \eqref{eq:state-combined} for almost every $t\in(0,T)$. Possible atoms at the terminal time $T$ do not enter inequalities asserted for almost every $t<T$ and are not detected by test functions supported in $[0,T)$. They may therefore be retained in the measures $\mu_{\rm pl}$ and $\mu_{\rm th}$ without changing the formulation.

\noindent
All requirements of Definition~\ref{def:state} have now been established. It is important to retain the distinction between the two thermal quantities appearing in the limit: the positive temperature is $\theta=e^\tau$, whereas the thermal energy entering the momentum and energy balances is represented by the measure $\Theta$. In general, $\Theta(t)=e^{\tau(t)}\dd x+\lambda_t$ and the possible excess measure $\lambda_t$ records thermal concentration not seen by the almost-everywhere limit of the temperatures. The entropy--energy formulation preserves this distinction throughout.

\section{Weighted plastic flow and lifted maximum plastic work}
\label{sec:plastic}
The construction of the relaxed entropy--energy state is now complete. We next establish the additional weighted and lifted constitutive structure generated by the approximation. The question is whether the exact normality encoded on the lifted fibres descends to a macroscopic Prandtl--Reuss relation between the limiting stress $\TT$ and the distributional plastic balance $\Xi$. The results below establish weighted compactness, lifted normality, and maximum plastic work; the limitations of macroscopic identification are discussed in Section~\ref{sec:normality}.

\noindent
Set $\omega_n=(\dev\TT^n-Z_n)/d_n$ where $d_n>0$, and choose any fixed unit deviator where $d_n=0$. Every measure lifted below vanishes where $\xi^n=0$, so that this choice has no effect. Lemma~\ref{lem:geometry} gives
\begin{equation}\label{eq:exact-projected-work}
\xi^n\in N_{r_nC_0}(Z_n)\,,\qquad \ip{Z_n}{\xi^n}=r_nh_0(\omega_n)|\xi^n|\,, \qquad q^n=r_nh_0(\omega_n)|\xi^n|+d_n|\xi^n|\,.
\end{equation}
The first term is the plastic work associated with the projected stress \(Z_n\in r_nC_0\) and therefore retains the exact support-function and normality structure of the yield set. The second term measures the additional work generated by the overstress \(d_n\); although \(d_n\to0\) strongly, this contribution cannot a priori be discarded because \(|\xi^n|\) is not uniformly integrable.

\subsection{Weighted compactness and support-function duality}
\begin{lemma}[Support-weighted plastic estimate]\label{lem:weighted}
The measures
\begin{equation}\label{eq:weighted-approx}
\Lambda_n:=r_n\xi^n\,\dd x\dd t\in\M(\qbar;\Sdev)\,,\qquad \nu^n:=r_n|\xi^n|\,\dd x\dd t\in\M^+(\qbar)\,,
\end{equation}
and $\chi_n:=r_nh_0(\omega_n)|\xi^n|\,\dd x\dd t\in\M^+(\qbar)$ are uniformly bounded. More precisely,
\begin{equation}\label{eq:weighted-pointwise}
r_n|\xi^n| \le r_nh_0(\omega_n)|\xi^n| \le \theta_c\frac{q^n}{\theta^n} \qquad\text{a.e. in }Q\,.
\end{equation}
Moreover, for $A_{n,\delta}:=\{(t,x)\in Q:r_n(t,x)\ge\delta\}$ $(0<\delta<r_{\max})$ one has
\begin{equation}\label{eq:subcritical-mass}
\int_{A_{n,\delta}}|\xi^n|\dd x\dd t \le \frac{C}{\delta}\,.
\end{equation}
\end{lemma}
\begin{proof}
On the set $\{r_n>0\}$, assumption~\eqref{eq:radius} implies $\theta^n<\theta_c$. Moreover, by \eqref{eq:exact-projected-work}, $r_nh_0(\omega_n)|\xi^n| = \ip{Z_n}{\xi^n} \le q^n$. Hence
\[
r_nh_0(\omega_n)|\xi^n| \le q^n \le \theta_c\frac{q^n}{\theta^n} \qquad\text{on }\{r_n>0\}\,.
\]
On $\{r_n=0\}$ the first two quantities in \eqref{eq:weighted-pointwise} vanish. Since $h_0(\omega_n)\ge1$, we also have $r_n|\xi^n| \le r_nh_0(\omega_n)|\xi^n|$. The uniform bound on $q^n/\theta^n$ from Lemma~\ref{lem:uniform} therefore yields the uniform measure bounds. Finally, on $A_{n,\delta}$ we have $r_n\ge\delta$, and thus
\[
\delta\int_{A_{n,\delta}}|\xi^n|\,\dd x\dd t \le \int_{A_{n,\delta}}r_n|\xi^n|\,\dd x\dd t
\le C\,,
\]
which proves \eqref{eq:subcritical-mass}.
\end{proof}

\noindent
After passing to a subsequence, we may therefore assume that
\[
\Lambda_n\weakstar\Lambda \quad\text{in }\M(\qbar;\Sdev)\,,\qquad \nu^n\weakstar\nu\,,\qquad \chi_n\weakstar\chi \quad\text{in }\M(\qbar)\,.
\]
Lower semicontinuity of the total variation and passage to the limit in the nonnegative measure inequalities in \eqref{eq:weighted-pointwise} give
\begin{equation}\label{eq:weighted-dominations}
|\Lambda| \le \nu \le \chi \le \theta_c\mu_{\rm pl}\,.
\end{equation}
The estimate $|\Lambda|\le\chi$ controls only the Euclidean total variation and does not yet retain the geometry of the yield body $C_0$. To recover this anisotropic information, we use the following dual characterization of the support-function variation, which allows it to be expressed through pairings with continuous fields taking values in $aC_0$. The measure-theoretic ingredients used in the proof are standard; see, for instance, \cite{AFP}.

\begin{lemma}[Support-function duality for measures]
\label{lem:measure-duality}
Let $P=\zeta\,|P|\in\M(X;\Sdev)$ be a finite $\Sdev$-valued Radon measure on a compact metric space $X$, and let
$\varphi\in C(X)$ be nonnegative. Then, for every $a\ge0$,
\begin{equation*}
a\int_X\varphi h_0(\zeta)\dd|P| = \sup\left\{ \int_X\varphi\ip{S}{\dd P}: S\in C(X;\Sdev)\,,\
S(y)\in aC_0\ \text{for every }y\in X \right\}.
\end{equation*}
If $X$ is an open cylinder, the same identity holds for $\varphi\in C_c(X)$ with competitors $S\in C_c(X;\Sdev)$ satisfying $S(y)\in aC_0$.
\end{lemma}

\begin{proof}
For every admissible field $S$ we have $\ip{S}{\zeta} \le a h_0(\zeta)$ $(|P|\text{-a.e.})$ and hence
\[
\int_X\varphi\ip{S}{\dd P}= \int_X\varphi\,\ip{S}{\zeta}\dd|P| \le a\int_X\varphi\,h_0(\zeta)\dd|P|\,.
\]
This proves one inequality. For the converse, choose a countable dense family $(S_j)_{j\in\mathbb N}$ in $aC_0$. Then $\sup_{j\in\mathbb N}\ip{S_j}{\zeta} = a h_0(\zeta)$ and finite maxima increase pointwise to the right-hand side. A measurable choice of maximizing indices therefore gives simple $aC_0$-valued fields whose pairings with $\zeta$ converge increasingly to $a h_0(\zeta)$. By regularity of the finite measure $\varphi|P|$ and Lusin approximation, these simple fields can be approximated in $L^1(\varphi|P|)$ by continuous fields. Projecting the approximants onto the closed convex set $aC_0$ preserves this convergence, since the metric projection is nonexpansive and fixes every $aC_0$-valued field. Passing to the limit gives the reverse inequality.

If $X$ is an open cylinder and $\varphi\in C_c(X)$, one additionally
multiplies the continuous approximants by a compactly supported cutoff
which equals one on $\supp\varphi$. Since $0\in aC_0$ and $aC_0$ is convex,
this operation preserves the pointwise constraint. The case $a=0$ is
immediate.
\end{proof}
\noindent
Applied to \(P=\Lambda\), this formula allows the pointwise \(C_0\)-constraint on the approximate weighted flows to pass to the limit and yields the anisotropic domination \(|\Lambda|_{C_0}\le\chi\).

\subsection{Weighted Minty inequality}

The homothetic structure of the yield sets provides comparison points even at complete collapse. Indeed, if $\bar Z\in C_0$ and $\eta\in N_{C_0}(\bar Z)$, then
\begin{equation}\label{eq:scaled-normal}
\eta\in N_{\rho C_0}(\rho\bar Z) \qquad\text{for every }\rho\ge0\,.
\end{equation}
For $\rho>0$ this follows directly from the definition of the normal cone, while for $\rho=0$ one has $N_{\{0\}}(0)=\Sdev$. Thus the same comparison remains admissible when the yield section collapses. This observation underlies the following limiting Minty inequality and is analogous to the homothetic comparisons used for temperature-dependent yield sets in \cite{BartelsRoubicek2008}.

\begin{proposition}[Weighted Minty inequality for the limiting measures]
\label{prop:minty}
Let $\varphi\in C(\qbar)$ be nonnegative, and let $\bar Z,\eta\in C(\qbar;\Sdev)$ satisfy $\bar Z(y)\in C_0$ and $\eta(y)\in N_{C_0}(\bar Z(y))$ for every $y\in\qbar$. Then
\begin{equation}\label{eq:weighted-minty}
\int_{\qbar}\varphi\dd\chi - \int_{\qbar}\varphi\ip{\bar Z}{\dd\Lambda} - \int_Q \varphi\ip{\eta}{\dev\TT-r(e^\tau)\bar Z}\dd x\dd t \ge0\,.
\end{equation}
Moreover,
\begin{equation}\label{eq:weighted-anisotropic}
|\Lambda|_{C_0}\le\chi \qquad\text{as measures on }\qbar\,.
\end{equation}
\end{proposition}
\begin{proof}
For each $n$, apply \eqref{eq:scaled-normal} with $\rho=r_n$. Since $\xi^n\in N_{r_nC_0}(Z_n)$ and $\eta\in N_{r_nC_0}(r_n\bar Z)$, then the monotonicity of the normal-cone operator gives
\[
\ip{\xi^n-\eta}{Z_n-r_n\bar Z}\ge0 \qquad\text{a.e. in }Q\,.
\]
Multiplying by $\varphi\ge0$ and integrating over $Q$, we obtain
\[
\begin{aligned}
0\le{} \int_Q\varphi\ip{\xi^n}{Z_n}\dd x\dd t - 
\int_Q\varphi\,r_n\ip{\xi^n}{\bar Z}\dd x\dd t - \int_Q\varphi\ip{\eta}{Z_n}\dd x\dd t +
\int_Q\varphi\, r_n\ip{\eta}{\bar Z} \dd x\dd t\,.
\end{aligned}
\]
By the definitions of $\chi_n$ and $\Lambda_n$, this is equivalently
\[
0\le \int_{\qbar}\varphi\dd\chi_n - \int_{\qbar}\varphi\ip{\bar Z}{\dd\Lambda_n} -
\int_Q\varphi\ip{\eta}{Z_n}\dd x\dd t + \int_Q\varphi\,r_n\ip{\eta}{\bar Z}\dd x\dd t\,.
\]
The first two terms pass to the limit by the weak-star convergence of $\chi_n$ and $\Lambda_n$. Since $Z_n\weakstar\dev\TT$ in $L^\infty(Q;\Sdev)$, the third term converges to $-\int_Q\varphi\,\ip{\eta}{\dev\TT}\,\dd x\dd t$. The strong convergence $r_n\longrightarrow r(e^\tau)$ in every finite $L^p(Q)$ yields convergence of the last term. Hence
\eqref{eq:weighted-minty} follows. 

\noindent
To obtain the anisotropic estimate, choose $\eta=0$. Since $0\in N_{C_0}(\bar Z)$ for every $\bar Z\in C_0$, \eqref{eq:weighted-minty} gives
\[
\int_{\qbar}\varphi\ip{\bar Z}{\dd\Lambda} \le \int_{\qbar}\varphi\dd\chi
\]
for every nonnegative $\varphi\in C(\qbar)$ and every continuous field $\bar Z$ with values in $C_0$. Applying Lemma~\ref{lem:measure-duality} with $P=\Lambda$ and $a=1$, and taking the supremum over all such fields $\bar Z$, gives $\int_{\qbar}\varphi\,\dd|\Lambda|_{C_0} \le \int_{\qbar}\varphi\,\dd\chi$. Since this holds for every nonnegative $\varphi\in C(\qbar)$, we obtain \eqref{eq:weighted-anisotropic}.
\end{proof}

\noindent
Under the stress admissibility $\dev\TT\in r(e^\tau)C_0$ the weighted Minty inequality is in fact equivalent to the anisotropic support-function bound \eqref{eq:weighted-anisotropic}. Indeed, if
$\eta\in N_{C_0}(\bar Z)$, then
\[
\ip{\eta}{\dev\TT-r(e^\tau)\bar Z}\le0\,,
\]
while \eqref{eq:weighted-anisotropic} gives $\int_{\qbar}\varphi\,\ip{\bar Z}{\dd\Lambda} \le \int_{\qbar}\varphi\,\dd\chi$. Together these inequalities recover \eqref{eq:weighted-minty}; the converse follows from \eqref{eq:weighted-minty} by choosing $\eta=0$ and applying Lemma~\ref{lem:measure-duality}. Thus the monotonicity of the approximate normal-cone relations survives the limit in a weighted anisotropic form, which remains meaningful through complete collapse. It does not, however, recover a macroscopic normality relation between $\TT$ and $\Xi$, nor does it retain the full stress--direction correlation. The lifted representation constructed below records precisely this finer information before barycentric cancellation.

\subsection{Global lifts and entropy production}

The measures $\Lambda$, $\nu$, and $\chi$ retain the size and the support-function structure of the weighted plastic flow, but they no longer record the pointwise correlation between the projected stress $Z_n$ and the plastic direction $\omega_n$. To preserve this information, we lift the approximate measures to the product of spacetime and the compact fibre $\mathcal F$. For $y=(t,x)\in Q$, set $L_n(y)=\bigl(y,r_n(y),Z_n(y),\omega_n(y)\bigr)\in\qbar\times\mathcal F$ and define
\begin{equation}\label{eq:global-lift-approx}
\Gamma_{{\rm pl},n} =(L_n)_\#\left(\frac{q^n}{\theta^n}\dd x\dd t\right)\,, \qquad
\widetilde\Gamma_n=(L_n)_\#\nu^n\,.
\end{equation}
Thus $\Gamma_{{\rm pl},n}$ is weighted by the full plastic entropy production, whereas $\widetilde\Gamma_n$ is weighted by $r_n|\xi^n|$. The distinction is important when the yield radius collapses: the first measure may still carry mass on $\{r_n=0\}$, while the second one vanishes there.

\begin{proposition}[Global lifted maximum plastic work]\label{prop:global-lifts}
Along a further subsequence,
\[
\Gamma_{{\rm pl},n}\weakstar\Gamma_{\rm pl}\,, \qquad \widetilde\Gamma_n\weakstar\widetilde\Gamma
\quad\text{in }\M^+(\qbar\times\mathcal F)\,.
\]
Both limits are concentrated on $\qbar\times\mathcal G$ and satisfy all identities and the domination in \eqref{eq:main-lifts}.
\end{proposition}
 \begin{proof}
By Lemmas~\ref{lem:uniform} and \ref{lem:weighted}
\[
\Gamma_{{\rm pl},n}(\qbar\times\mathcal F) =\int_Q\frac{q^n}{\theta^n}\dd x\dd t\le C\,, \qquad
\widetilde\Gamma_n(\qbar\times\mathcal F) =\nu^n(\qbar)\le C\,.
\]
Since $\qbar\times\mathcal F$ is compact, weak-star compactness of finite Radon measures gives the asserted convergences after extraction. Whenever one of the weights in \eqref{eq:global-lift-approx} is nonzero, one has $\xi^n\ne0$. By the definition of the Norton--Hoff flow, this implies $d_n>0$, and hence $\omega_n=\frac{\xi^n}{|\xi^n|}$. Together with \eqref{eq:exact-projected-work}, this gives $Z_n\in r_nC_0$ and $\ip{Z_n}{\omega_n}=r_nh_0(\omega_n)$. Consequently both approximate lifts are concentrated on $\qbar\times\mathcal G$. The set $\mathcal G$ is closed. Indeed, if
\[
(\rho_j,Z_j,\omega_j)\longrightarrow(\rho,Z,\omega) \quad\text{in }\mathcal F\,, \qquad
(\rho_j,Z_j,\omega_j)\in\mathcal G\,,
\]
then, by compactness of $C_0$, the relations $Z_j\in\rho_jC_0$ imply $Z\in\rho C_0$. Moreover, continuity of the scalar product and of the support function $h_0$ yields $\ip{Z}{\omega}=\rho h_0(\omega)$. Thus $\mathcal G$ is closed. Since both approximate lifts are concentrated on $\qbar\times\mathcal G$, their weak-star limits are concentrated there as well. Consequently, $\ip{Z}{\omega} =\rho h_0(\omega) =\max_{S\in\rho C_0}\ip{S}{\omega}$ holds $\Gamma_{\rm pl}$-a.e. and $\widetilde\Gamma$-a.e.

\noindent
The identities for the marginals follow directly from the definition of the push-forward measures. Indeed, $\pi_\#\Gamma_{{\rm pl},n} =\frac{q^n}{\theta^n}\dd x\dd t$, $\pi_\#\widetilde\Gamma_n=\nu^n$ and therefore
\[
\pi_\#\Gamma_{\rm pl}=\mu_{\rm pl}\,, \qquad \pi_\#\widetilde\Gamma=\nu.
\]
Similarly, testing the second lift with the continuous fibre functions $\omega$ and $h_0(\omega)$ gives $\pi_\#(\omega\widetilde\Gamma_n)=\Lambda_n$ and $\pi_\#(h_0(\omega)\widetilde\Gamma_n)=\chi_n$. Passing to the limit yields the corresponding identities for $\Lambda$ and $\chi$. It remains to compare the two lifts. For every nonnegative
$\Psi\in C(\qbar\times\mathcal F)$, \eqref{eq:weighted-pointwise} gives
\[
\begin{aligned}
\int_{\qbar\times\mathcal F} \Psi h_0(\omega)\dd\widetilde\Gamma_n &= \int_Q \Psi(L_n(y))\,r_nh_0(\omega_n)|\xi^n| \dd x\dd t  \le \theta_c \int_Q \Psi(L_n(y))\,\frac{q^n}{\theta^n} \dd x\dd t\\
&= \theta_c \int_{\qbar\times\mathcal F} \Psi\dd\Gamma_{{\rm pl},n}\,.
\end{aligned}
\]
Passing to the limit proves $h_0(\omega)\widetilde\Gamma\le\theta_c\Gamma_{\rm pl}$.
\end{proof}

\noindent
The decomposition in \eqref{eq:exact-projected-work} also separates the plastic entropy production into a part generated by the projected stress and a part generated by the overstress. More precisely, set
\[
\mu_{{\rm proj},n} = \frac{r_nh_0(\omega_n)|\xi^n|}{\theta^n} \dd x\dd t\,, \qquad
\mu_{{\rm over},n} = \frac{d_n|\xi^n|}{\theta^n} \dd x\dd t\,.
\]
Then
\[
\frac{q^n}{\theta^n}\dd x\dd t = \mu_{{\rm proj},n}+\mu_{{\rm over},n}\,.
\]

\begin{proposition}[Projected and overstress entropy production]\label{prop:production-split}
Along a further subsequence there exist nonnegative measures $\mu_{\rm proj},\mu_{\rm over}\in\M^+(\qbar)$ such that
\begin{equation}\label{eq:production-split}
\begin{aligned}
\mu_{{\rm proj},n}&\weakstar\mu_{\rm proj}\,,\qquad \mu_{{\rm over},n}&\weakstar\mu_{\rm over}\,,
\qquad \mu_{\rm pl}=\mu_{\rm proj}+\mu_{\rm over}\,.
\end{aligned}
\end{equation}
Moreover, $\chi\le\theta_c\mu_{\rm proj}$.
\end{proposition}

\begin{proof}
Both approximate measures are nonnegative and their sum is
$q^n/\theta^n\,\dd x\dd t$. The entropy estimate therefore gives uniform bounds on their masses, and weak-star compactness yields the first two convergences in \eqref{eq:production-split}. Passing to the limit in their sum gives $\mu_{\rm pl}=\mu_{\rm proj}+\mu_{\rm over}$. Since
\[
\chi_n
=
r_nh_0(\omega_n)|\xi^n|\dd x\dd t
=
\theta^n\mu_{{\rm proj},n}
\]
and $r_n>0$ implies $\theta^n<\theta_c$, while both sides vanish when $r_n=0$, one has $\chi_n\le\theta_c\mu_{{\rm proj},n}$. The asserted inequality follows by passage to the limit.
\end{proof}

\noindent
The measure $\mu_{\rm over}$ cannot in general be discarded solely from the strong Lebesgue convergence $d_n\to0$. Indeed, the density defining $\mu_{{\rm over},n}$ contains the factor $|\xi^n|/\theta^n$, for which no corresponding equiintegrability is available. Thus a quantity that vanishes strongly with respect to Lebesgue measure may still contribute after multiplication by the concentrating plastic rate. In particular, $\chi$, $\mu_{\rm proj}$, and $\mu_{\rm pl}$ represent different limiting quantities.

\noindent
The information retained by $\widetilde\Gamma$ can also be read directly through its disintegration with respect to its base marginal $\nu$. Write $\widetilde\Gamma (\dd y,\dd\rho,\dd Z,\dd\omega) = \widetilde\gamma_y(\dd\rho,\dd Z,\dd\omega)\,\nu(\dd y)$, where $\widetilde\gamma_y$ is a probability measure for $\nu$-almost every $y$. Then
\begin{equation}\label{eq:global-disintegration}
\Lambda = \left(\int\omega\dd\widetilde\gamma_y\right)\nu\,, \qquad \chi = \left(\int h_0(\omega)\dd\widetilde\gamma_y\right)\nu\,.
\end{equation}
Thus the density of $\Lambda$ with respect to $\nu$ is the barycentre of the plastic directions, whereas the density of $\chi$ is the corresponding average support value. By convexity and positive homogeneity of $h_0$, $h_0\left(\int\omega\dd\widetilde\gamma_y\right) \le \int h_0(\omega)\dd\widetilde\gamma_y$, which recovers \eqref{eq:weighted-anisotropic}. Likewise, $\left|\int\omega\dd\widetilde\gamma_y\right|\le1$ gives $|\Lambda|\le\nu$, while $h_0(\omega)\ge1$ gives $\nu\le\chi$.

\noindent
The global lift therefore preserves the exact maximum-plastic-work relation at the fibre level, while the measures $\Lambda$, $\nu$, and $\chi$ retain only its barycentric consequences. To recover finite unweighted plastic flow, one must restrict to regions where the yield radius stays uniformly away from zero; this is the purpose of the next subsection.

\subsection{Plastic flow on nondegenerate approximate regions}

The weighted estimate of Lemma~\ref{lem:weighted} does not in general provide a bound for the full plastic rate $\xi^n$, because the radius $r_n$ may approach zero. Such a bound is recovered, however, on every region where the approximate yield radius stays uniformly positive. For $\delta\in(0,r_{\max})$, recall $A_{n,\delta}=\{(t,x)\in Q:r_n(t,x)\ge\delta\}$. Then
\begin{equation}\label{eq:subcritical-unweighted}
\int_{A_{n,\delta}}|\xi^n|\dd x\dd t \le \frac{1}{\delta}\int_Q r_n|\xi^n|\dd x\dd t \le \frac{C}{\delta}\,.
\end{equation}
Thus, away from the collapsing part of the yield set, the unweighted plastic flow is bounded in the space of finite measures. We therefore define
\begin{equation}\label{eq:delta-approx}
\begin{aligned}
m_{\delta,n} &=\mathbf1_{A_{n,\delta}}|\xi^n|\dd x\dd t\,,
& P_{\delta,n} &=\mathbf1_{A_{n,\delta}}\xi^n\dd x\dd t\,, \\
\Gamma_{\delta,n} &=(L_n)_\#m_{\delta,n}\,, & \mu_{\delta,n} &=\mathbf1_{A_{n,\delta}} r_nh_0(\omega_n)|\xi^n|\dd x\dd t\, .
\end{aligned}
\end{equation}
Here $m_{\delta,n}$ is the total variation measure of the plastic flow on the nondegenerate region, while $P_{\delta,n}$ retains its tensorial direction. The measure $\mu_{\delta,n}$ records the corresponding projected plastic work, and $\Gamma_{\delta,n}$ keeps, in addition, the fibre variables on which the exact normality relation is available. By Lemma~\ref{lem:weighted}, $m_{\delta,n}(\qbar)\le\frac{C}{\delta}$ and $\mu_{\delta,n}(\qbar)\le C$. Since the fibre $\mathcal F$ is compact, these bounds also yield compactness of the associated lifts.
\begin{proposition}[Level-restricted lifted normality and barycentres]
\label{prop:delta-lift}
For every fixed $\delta\in(0,r_{\max})$, after passing to a subsequence,
\[
m_{\delta,n}\weakstar\gamma_\delta\,,\qquad \Gamma_{\delta,n}\weakstar\Gamma_\delta\,,\qquad P_{\delta,n}\weakstar P_\delta\,,\qquad \mu_{\delta,n}\weakstar\mu_\delta\, .
\]
The measure $\Gamma_\delta$ is concentrated on $\{(y,\rho,Z,\omega): y\in\qbar,\ (\rho,Z,\omega)\in\mathcal G,\ \rho\ge\delta\}$. Moreover, $\gamma_\delta=\pi_\#\Gamma_\delta$ and
\begin{equation}\label{eq:delta-barycentres}
P_\delta=\pi_\#(\omega\Gamma_\delta)\,, \qquad \mu_\delta=\pi_\#(\rho h_0(\omega)\Gamma_\delta)\,.
\end{equation}
If $\eta_\delta:=\pi_\#(h_0(\omega)\Gamma_\delta)$, then
\begin{equation}\label{eq:delta-chain}
|P_\delta|_{C_0}\le\eta_\delta\,, \qquad \delta\eta_\delta\le\mu_\delta\le\chi\,.
\end{equation}
In particular, $\delta|P_\delta|_{C_0}\le\mu_\delta$. More precisely,
\begin{equation}\label{eq:delta-defects}
\mu_\delta-\delta|P_\delta|_{C_0} = (\mu_\delta-\delta\eta_\delta) +\delta(\eta_\delta-|P_\delta|_{C_0})\,,
\end{equation}
where both terms on the right-hand side are nonnegative measures. Finally, for every continuous field
$S:\qbar\to\delta C_0$ and every nonnegative
$\varphi\in C(\qbar)$,
\begin{equation}\label{eq:delta-variational}
\int_{\qbar}\varphi\,\ip{S}{\dd P_\delta} \le \int_{\qbar}\varphi\,\dd\mu_\delta\,.
\end{equation}
\end{proposition}

\begin{proof}
The uniform mass bounds above and compactness of the fibre yield simultaneous weak-star compactness of the four measures in \eqref{eq:delta-approx}. Since $\Gamma_{\delta,n}$ is concentrated on $\{(y,\rho,Z,\omega): (\rho,Z,\omega)\in\mathcal G,\ \rho\ge\delta\}$, which is closed, the same support property passes to $\Gamma_\delta$. The marginal and barycentre identities follow directly from the definition of the push-forward. Indeed, testing $\Gamma_{\delta,n}$ with functions of the form $\varphi(y)$, $\varphi(y)\omega$, $\varphi(y)\rho h_0(\omega)$ and passing to the limit gives respectively $\gamma_\delta=\pi_\#\Gamma_\delta$ and \eqref{eq:delta-barycentres}. To compare the barycentre $P_\delta$ with the positive measure $\eta_\delta$, disintegrate
\[
\Gamma_\delta(\dd y,\dd\rho,\dd Z,\dd\omega) = \gamma_{\delta,y}(\dd\rho,\dd Z,\dd\omega)\,
\gamma_\delta(\dd y)\,,
\]
where $\gamma_{\delta,y}$ is a probability measure for $\gamma_\delta$-almost every $y$. Set
\[
p_\delta(y)=\int\omega\,\dd\gamma_{\delta,y}\,.
\]
Then $P_\delta=p_\delta\,\gamma_\delta$ and, by positive homogeneity of $h_0$, $|P_\delta|_{C_0} = h_0(p_\delta)\,\gamma_\delta$. Since $h_0$ is convex, Jensen's inequality gives
\[
h_0(p_\delta) \le \int h_0(\omega)\,\dd\gamma_{\delta,y}\,,
\]
and therefore $|P_\delta|_{C_0}\le\eta_\delta$. On the support of $\Gamma_\delta$ one has $\rho\ge\delta$, and hence
\[
\mu_\delta-\delta\eta_\delta = \pi_\# \bigl((\rho-\delta)h_0(\omega)\Gamma_\delta\bigr) \ge0\,.
\]
This proves $\delta\eta_\delta\le\mu_\delta$. Furthermore, $\mu_{\delta,n}\le\chi_n$, and consequently $\mu_\delta\le\chi$. Combining these inequalities gives \eqref{eq:delta-chain}. Formula~\eqref{eq:delta-defects} is then obtained by adding and subtracting $\delta\eta_\delta$. The two terms have different origins. The first measures the excess of the actual fibre radius $\rho$ above the reference level $\delta$, whereas the second measures the loss produced by averaging different flow directions in the barycentre $P_\delta$. It remains to prove the variational inequality. On $A_{n,\delta}$ we have $r_n\ge\delta$. Since $0\in C_0$ and $C_0$ is convex, $\delta C_0\subset r_nC_0$. Hence every $S(y)\in\delta C_0$ is an admissible comparison stress for the normality relation $\xi^n\in N_{r_nC_0}(Z_n)$ and therefore
\[
\xi^n:S\le\xi^n:Z_n = r_nh_0(\omega_n)|\xi^n| \qquad\text{on }A_{n,\delta}\,.
\]
Multiplying by $\varphi\ge0$, integrating, and passing to the limit yields \eqref{eq:delta-variational}. Lemma~\ref{lem:measure-duality} then gives again $\delta|P_\delta|_{C_0}\le\mu_\delta$.
\end{proof}

\noindent
The preceding proposition describes the plastic flow on the moving sets $A_{n,\delta}=\{r_n\ge\delta\}$ and retains its normality structure through the lifted measure $\Gamma_\delta$. A stronger conclusion is available on a fixed open region which remains uniformly nondegenerate along the approximation. There the cutoff is eventually inactive, and the unweighted plastic flow can be identified directly with the distributional constitutive flow $\Xi$.

\begin{proposition}[Identification on a fixed nondegenerate region]
\label{prop:fixed-region-identification}
Let the approximating sequence be taken along the subsequence fixed in Theorem~\ref{thm:main}. Suppose that $U\subset Q$ is open and that, for some $\delta\in(0,r_{\max})$ and $n_0\in\mathbb N$,
\begin{equation}\label{eq:fixed-region-positive}
r_n\ge\delta \quad\text{a.e. in }U \qquad\text{for every }n\ge n_0\,.
\end{equation}
Then
\begin{equation}\label{eq:fixed-region-identification}
\bigl(\xi^n\dd x\dd t\bigr)\restriction_U \weakstar P_U \quad\text{in }\M(U;\Sdev)\,, \qquad
\Xi\restriction_U=P_U\,,
\end{equation}
along the already fixed subsequence and without any further extraction. Moreover,
\begin{equation}\label{eq:fixed-region-bounds}
|P_U|(U)\le\frac{C}{\delta}\,,
\qquad
\delta|P_U|_{C_0} \le\chi\restriction_U \le\theta_c\mu_{\rm pl}\restriction_U\,.
\end{equation}
If $0<\eta\le\delta$ and the level-restricted measure $P_\eta$ in \eqref{eq:main-local} has been selected, then
\begin{equation}\label{eq:fixed-region-level}
P_\eta\restriction_U = P_U = \Xi\restriction_U\,.
\end{equation}
\end{proposition}

\begin{proof}
By Lemma~\ref{lem:weighted} and \eqref{eq:fixed-region-positive},
\[
\sup_{n\ge n_0} \int_U|\xi^n|\dd x\dd t \le \frac1\delta \sup_{n\ge n_0} \int_Qr_n|\xi^n|\dd x\dd t \le\frac{C}{\delta}\,.
\]
Thus $\bigl(\xi^n\dd x\dd t\bigr)\restriction_U$ is weak-star relatively compact in $\M(U;\Sdev)$. Let $P_U$ be the limit of an arbitrary convergent subsequence. For every $\Phi\in C_c^\infty(U;\Sdev)$, the distributional convergence obtained in Section~\ref{sec:limit} gives
\[
\int_U\ip{\Phi}{\dd P_U} = \lim_n\int_U\ip{\Phi}{\xi^n}\dd x\dd t=\langle\Xi,\Phi\rangle\,.
\]
Consequently, $P_U=\Xi\restriction_U$ as distributions on $U$. In particular, the restriction of $\Xi$ to $U$ is represented by a finite Radon measure. Since every subsequential measure limit has the same value, uniqueness of the cluster point together with the uniform mass bound gives the convergence in \eqref{eq:fixed-region-identification} along the full subsequence already fixed in Theorem~\ref{thm:main}. The estimate $|P_U|(U)\le C/\delta$ follows at the same time.

To prove the anisotropic estimate, let
$\varphi\in C_c(U)$ be nonnegative. By weak-star lower semicontinuity and Lemma~\ref{lem:measure-duality},
\[
\delta\int_U\varphi\,\dd|P_U|_{C_0} \le \liminf_n \delta\int_U\varphi\,h_0(\xi^n)\dd x\dd t\,.
\]
On $\{\xi^n\ne0\}$ the direction of $\xi^n$ is $\omega_n$, and positive homogeneity of $h_0$ gives $h_0(\xi^n)=h_0(\omega_n)|\xi^n|$. Using $r_n\ge\delta$ on $U$, we obtain
\[
\begin{aligned}
\delta\int_U\varphi\,\dd|P_U|_{C_0} \le \liminf_n \int_U\varphi\,r_nh_0(\omega_n)|\xi^n|\dd x\dd t = \int_U\varphi\,\dd\chi.
\end{aligned}
\]
Hence $\delta|P_U|_{C_0}\le\chi\restriction_U$. The remaining inequality in \eqref{eq:fixed-region-bounds} follows from \eqref{eq:weighted-dominations}. Finally, if $0<\eta\le\delta$, then $\mathbf1_{\{r_n\ge\eta\}}\xi^n=\xi^n$ a.e. in $U$, for every $n\ge n_0$. Testing the convergence defining $P_\eta$ with fields compactly supported in $U$ therefore gives \eqref{eq:fixed-region-level}.
\end{proof}

\begin{remark}[Approximate and limiting nondegeneracy]
\label{rem:fixed-region-scope}
Condition~\eqref{eq:fixed-region-positive} is a sufficient local condition on the approximate radii and is not part of the global construction. It is stronger than the corresponding positivity condition for the limiting radius. In particular, $r(e^\tau)\ge\delta$ a.e. in $U$, together with the almost-everywhere and strong finite-$L^p$ convergence $r_n\to r(e^\tau)$ does not imply
\eqref{eq:fixed-region-positive}. Therefore Proposition~\ref{prop:fixed-region-identification} does not identify the general level-restricted measure $P_\delta$ with $\Xi\restriction_{\{r(e^\tau)\ge\delta\}}$. Even on a region where $\Xi$ is identified as a finite measure, the existence of a compatible pairing between $\dev\TT$ and $\Xi$ and the corresponding macroscopic maximum-plastic-work identity remain separate questions.
\end{remark}

\noindent
For the von Mises ball, $h_0(V)=|V|$ and $|P_\delta|_{C_0}=|P_\delta|$. Hence \eqref{eq:delta-chain} reduces to $\delta|P_\delta|\le\mu_\delta$. For a general convex body $C_0$, the measure $|P_\delta|_{C_0}$ retains the anisotropic geometry of the yield surface.

\noindent
The construction above can be carried out simultaneously on any prescribed countable family of positive levels by a diagonal argument. The next subsection removes this dependence on separate levelwise extractions and organizes the restrictions through a single locally finite lift on the positive-radius fibres.

\subsection{A canonical lift on the positive-radius fibres}

The preceding construction was performed separately on the sets $\{r_n\ge\delta\}$ for fixed positive levels $\delta$. Although this recovers the unweighted plastic flow away from the collapsing regime, a separate extraction may in principle depend on the chosen level. We now show that all positive levels can instead be organized through the single support-weighted lift $\widetilde\Gamma$.

\noindent
Set
\[
X=\qbar\times\mathcal F\,, \qquad X_+=X\cap\{\rho>0\}\,.
\]
Recall that the fibre variable $\rho$ corresponds, at the approximate level, to the yield radius $r_n$. Since $\widetilde\Gamma_n$ carries the weight $r_n$, division by $\rho$ removes this weight whenever $\rho>0$. This motivates
\[
\Gamma_n^+ := \rho^{-1}\widetilde\Gamma_n\restriction_{X_+}\,, \qquad \Gamma^+ := \rho^{-1}\widetilde\Gamma\restriction_{X_+}\,.
\]
The factor $\rho^{-1}$ may become singular as $\rho\downarrow0$, so $\Gamma^+$ need not be finite on all of $X_+$. It is, however, finite on every compact subset of $X_+$. We denote the corresponding class by $\M_{\rm loc}^+(X_+)$.

\begin{proposition}[Canonical positive-radius lift for a fixed subsequence]
\label{prop:positive-lift}
Let $\widetilde\Gamma$ be the limit obtained in Proposition~\ref{prop:global-lifts}. Then $\Gamma^+ = \rho^{-1}\widetilde\Gamma\restriction_{X_+} \in\M_{\rm loc}^+(X_+)$ and $\Gamma^+$ is $\sigma$-finite. Moreover, along the same subsequence,
\begin{equation}\label{eq:positive-vague}
\Gamma_n^+ \longrightarrow\Gamma^+ \qquad \text{against every }\Psi\in C_c(X_+)\,.
\end{equation}
Define the set of exceptional levels by $\mathcal E = \bigl\{ \delta\in(0,r_{\max}): \widetilde\Gamma(\{\rho=\delta\})>0 \bigr\}$. Then $\mathcal E$ is at most countable. For every
$\delta\in(0,r_{\max})\setminus\mathcal E$, no further extraction is required and
\begin{equation}\label{eq:positive-level-limit}
\Gamma_{\delta,n}\weakstar \Gamma_\delta = \Gamma^+\restriction_{\{\rho\ge\delta\}} \qquad\text{in }\M^+(X)\,.
\end{equation}
In particular, all barycentric identities and estimates of Proposition~\ref{prop:delta-lift} hold for these limits. If $\delta_j\downarrow0$ through nonexceptional levels, then $\Gamma_{\delta_j}\uparrow\Gamma^+$ as positive measures on $X_+$.
\end{proposition}

\begin{proof}
We first prove local finiteness. Let $K\subset X_+$ be compact. Since $K$ does not meet $\{\rho=0\}$, $\delta_K:=\min_K\rho>0$. Consequently,
\[
\Gamma^+(K) = \int_K\frac1\rho\,\dd\widetilde\Gamma \le \frac1{\delta_K}\widetilde\Gamma(K) <\infty\,.
\]
Thus $\Gamma^+\in\M_{\rm loc}^+(X_+)$. Moreover,
\[
X_+= \bigcup_{j=1}^\infty \left(X\cap\left\{\rho\ge\frac{r_{\max}}{j}\right\}\right)\,,
\]
and $\Gamma^+$ has finite mass on each set in this union. Hence $\Gamma^+$ is $\sigma$-finite. Since $\widetilde\Gamma$ is concentrated on the normality graph $\mathcal G$, the same is true of $\Gamma^+$ on $X_+$.

\noindent
Let now $\Psi\in C_c(X_+)$. Its compact support stays a positive distance away from $\{\rho=0\}$. Therefore $\Psi/\rho$ is continuous on its support and extends by zero to a continuous function on $X$. Using $\widetilde\Gamma_n\weakstar\widetilde\Gamma$, we obtain
\[
\int_{X_+}\Psi\,\dd\Gamma_n^+ = \int_X\frac{\Psi}{\rho}\,\dd\widetilde\Gamma_n \longrightarrow \int_X\frac{\Psi}{\rho}\,\dd\widetilde\Gamma = \int_{X_+}\Psi\,\dd\Gamma^+.
\]
This proves \eqref{eq:positive-vague}. We next relate $\Gamma^+$ to the fixed-level lifts of the preceding subsection. At the approximate level,
\begin{equation}\label{eq:positive-cutoff-identity}
\Gamma_{\delta,n} = \frac{\mathbf1_{\{\rho\ge\delta\}}}{\rho}\, \widetilde\Gamma_n\,,
\end{equation}
where the multiplier is set equal to zero on $\{\rho<\delta\}$. Thus $\Gamma_{\delta,n}$ is simply the positive-radius lift with all fibres below the level $\delta$ removed. The only discontinuity of the multiplier in \eqref{eq:positive-cutoff-identity} occurs at $\rho=\delta$. Hence, if $\widetilde\Gamma(\{\rho=\delta\})=0$, the cutoff can be approximated from above and below by continuous functions. Passing to the limit in \eqref{eq:positive-cutoff-identity} gives
\[
\Gamma_{\delta,n}\weakstar \frac{\mathbf1_{\{\rho\ge\delta\}}}{\rho}\, \widetilde\Gamma = \Gamma^+\restriction_{\{\rho\ge\delta\}}\,,
\]
which proves \eqref{eq:positive-level-limit} without any further subsequence.

\noindent
The exceptional levels are precisely the atoms in $(0,r_{\max})$ of the finite measure $(\operatorname{pr}_\rho)_\#\widetilde\Gamma$. A finite measure has at most countably many atoms, so $\mathcal E$ is at most countable. For every nonexceptional level the barycentric identities and bounds now follow from Proposition~\ref{prop:delta-lift}. Finally, if
$\delta_j\downarrow0$, then $\{\rho\ge\delta_j\}\uparrow X_+$ and therefore $\Gamma_{\delta_j} = \Gamma^+\restriction_{\{\rho\ge\delta_j\}} \uparrow \Gamma^+$.
\end{proof}

\begin{remark}[Exceptional levels and scope of the positive-radius lift]
\label{rem:positive-scope}
At an exceptional level $\delta\in\mathcal E$, the limit of $\Gamma_{\delta,n}$ is not determined by $\widetilde\Gamma$ alone. Indeed, after a further extraction, every possible limit has the form
\[
\Gamma_\delta = \rho^{-1}\widetilde\Gamma\restriction_{\{\rho>\delta\}} + \delta^{-1}B_\delta\,,
\qquad 0\le B_\delta \le \widetilde\Gamma\restriction_{\{\rho=\delta\}}\,.
\]
The part strictly above $\delta$ is uniquely determined, while the measure $B_\delta$ records the portion of the limiting mass at the level $\rho=\delta$ which is selected by the approximate cutoffs. A diagonal argument can be used for any prescribed countable collection of exceptional levels. The measure $\Gamma^+$ is canonical once the limit $\widetilde\Gamma$ and hence the underlying subsequence have been fixed. This does not imply that it is independent of the original regularization or of the selected subsequence.

\noindent
Moreover, $\Gamma^+$ need not have finite total mass. Its projection onto the physical cylinder $\qbar$ need not be locally finite either, because arbitrarily large unweighted plastic mass may accumulate on fibres with $\rho\downarrow0$. The identity $\rho\Gamma^+ = \widetilde\Gamma\restriction_{X_+}$ also shows that $\Gamma^+$ contains no information about possible mass of $\widetilde\Gamma$ carried by the collapsed fibre $\{\rho=0\}$.

Finally, the fibre variables retained by $\Gamma^+$ should not be identified automatically with the macroscopic limiting fields. In particular, the lift construction alone does not identify $\rho$ with $r(e^\tau)$ with respect to the plastic measure, nor does it identify a barycentre of $\Gamma^+$ with the distributional constitutive flow $\Xi$. The identification $P_U=\Xi\restriction_U$ obtained in Proposition~\ref{prop:fixed-region-identification} uses,
in addition to measure compactness, the distributional constitutive balance and the uniform approximate lower bound \eqref{eq:fixed-region-positive}. The monotonicity $\Gamma_{\delta_j}\uparrow\Gamma^+$ concerns the positive lifted measures and does not imply monotonicity of the vector measures $P_{\delta_j}$.
\end{remark}

\begin{proof}[Completion of the proof of Theorem~\ref{thm:main}]
Proposition~\ref{prop:fixed-level} provides the approximate solutions. Lemma~\ref{lem:uniform}, Proposition~\ref{prop:thermal}, and the limit arguments of Section~\ref{sec:limit} yield the relaxed entropy--energy state and item~\textup{(i)}. Propositions~\ref{prop:overstress} and \ref{prop:admissibility} give item~\textup{(ii)}. Equation~\eqref{eq:Xi-bound}, Lemma~\ref{lem:weighted}, and Proposition~\ref{prop:minty} give item~\textup{(iii)}, while Propositions~\ref{prop:global-lifts} and \ref{prop:production-split} give item~\textup{(iv)}.

\noindent
Finally, Proposition~\ref{prop:delta-lift} provides the level-restricted lifts, and Proposition~\ref{prop:positive-lift} organizes all nonexceptional positive levels through the single locally finite measure $\Gamma^+$. A further diagonal extraction is needed only when a prescribed countable collection contains exceptional levels. This yields all the convergence, barycentric identities, and bounds stated in item~\textup{(v)} and completes the proof.

\noindent
Proposition~\ref{prop:fixed-region-identification} gives, in addition, the identification of the level-restricted flow with the full distributional constitutive flow on fixed regions which remain uniformly nondegenerate along the approximation. All extractions above are taken along subsequences of the previously fixed sequence and therefore preserve the earlier convergences.
\end{proof}

\begin{definition}[Approximation-generated state]
\label{def:generated}
A relaxed entropy--energy state in the sense of Definition~\ref{def:state}, together with its associated weighted and lifted constitutive measures, is called \emph{approximation-generated} if it is obtained from solutions of the approximate system \eqref{eq:approx}, with the regularized data constructed in Section~\ref{sec:approximation}, along one common subsequence. More precisely, along this subsequence the mechanical convergences \eqref{eq:mechanical-limits}, the thermal convergences \eqref{eq:main-tau}
\[
\theta^n\dd x\weakstar\Theta \qquad\text{in } L^\infty_{w*}(0,T;\M(\overline\Omega))\,,
\]
the production limits \eqref{eq:production-limits}, and all weighted and lifted limits constructed in Section~\ref{sec:plastic} are required to hold simultaneously.
\end{definition}

The limiting algebraic relations satisfied by the weighted and lifted measures do not by themselves imply macroscopic normality; see Remark~\ref{rem:abstract-lift-insufficient}. We neither prove uniqueness of the generated state nor independence of the regularization or of the subsequence used in its construction.

\section{Macroscopic versus lifted normality}
\label{sec:normality}

Section~\ref{sec:plastic} establishes exact normality and maximum plastic work on the lifted fibres. It does not, however, identify a corresponding macroscopic stress--plastic-flow pairing. The reason is that the limiting constitutive objects are defined with respect to different reference structures. The full plastic balance yields $\Xi\in H^{-1}(Q;\Sdev)$, whereas the finite measures $\Lambda$ and $P_\delta$ retain weighted or level-restricted parts of the approximate plastic flow, and their lifts encode exact stress--direction normality. The strong convergence of $r_n$ and of the projection error is with respect to Lebesgue measure, while the lifted measures are generated by plastic densities that may concentrate.Consequently, the available compactness does not identify the fibre variables $\rho$ and $Z$ with $r(e^\tau)$ and a representative of $\dev\TT$ on the plastic reference measures.

\noindent
In particular, Theorem~\ref{thm:main} does not assert $P_\delta=\Xi\restriction_{\{r(e^\tau)\ge\delta\}}$. Indeed, the restriction of a general $H^{-1}$ distribution to a measurable set is not canonically defined. Moreover, even when a plastic measure is available on a given region, an $L^\infty$ stress has no canonical representative on its singular part. Recovering a macroscopic maximum-work relation would therefore require two additional ingredients: identification of the measure-valued plastic flow with the full constitutive balance and a compatible stress--plastic-flow pairing. The weighted Minty inequality retains the corresponding support bound through the independent measure $\chi$, but does not identify $\chi$ with such a pairing. In the coercive setting of \cite{BartelsRoubicek2008}, additional rate estimates and convergence of the mechanical product provide precisely this missing identification. Stress--plastic-strain duality and the associated boundary formulation are likewise essential components of the classical theories \cite{AnzellottiLuckhaus1987,DalMasoDeSimoneMora2006,BabadjianMifsud2017}; the lifted identities obtained here do not by themselves replace them.

\subsection{The intrinsic balance class is strictly broader}

\begin{proposition}[A smooth state without macroscopic normality]
\label{prop:smooth-nonnormal}
The intrinsic formulation admits smooth states that remain strictly inside a nondegenerate yield set while violating the Prandtl--Reuss normality relation. In particular, relaxed entropy--energy states need not be unique even for smooth data.
\end{proposition}

\begin{proof}
Let $\C$ be the identity, $C_0=\overline B_1^{\Sdev}(0)$, $r(\theta)=(2-\theta)_+$ and $f=0$, and fix $E\in\Sdev$ with $|E|=1$ and $0<a_0<1/2$. For the data $u_0=u_1=0$, $\TT_0=a_0E$ and $\theta_0=1$, set
\begin{equation}\label{eq:smooth-nonnormal}
\begin{gathered}
a(t)=a_0e^{-t}\,,\qquad u=0\,,\qquad \TT=a(t)E\,,\qquad \Xi=a(t)E\,,\\
\theta(t)=1+\frac{a_0^2-a(t)^2}{2}\,,\qquad \tau=\log\theta\,,\qquad \Theta(t)=\theta(t)\dd x\,,\\
\mu_{\rm pl}=\frac{a(t)^2}{\theta(t)}\dd x\dd t\,,\qquad \mu_{\rm th}=0\,.
\end{gathered}
\end{equation}
Since the fields are spatially homogeneous, the momentum equation is immediate, while
\[
\TT_t+\Xi=0\,,\qquad \theta_t=a(t)^2=\ip{\Xi}{\TT}\,,\qquad \tau_t=\frac{a(t)^2}{\theta(t)}\,.
\]
Hence the constitutive and entropy balances hold, and $\frac12a(t)^2+\theta(t)=\frac12a_0^2+1$ gives the energy identity. The remaining regularity, initial conditions, thermal splitting, and volumetric identity follow directly from \eqref{eq:smooth-nonnormal}. On the other hand,
\[
|\dev\TT|=a(t)<\frac12\,,\qquad r(\theta(t))=1-\frac12\bigl(a_0^2-a(t)^2\bigr)>\frac78\,.
\]
Thus $\dev\TT$ lies strictly inside $r(\theta)C_0$, so $N_{K(\theta)}(\TT)=\{0\}$, whereas $\Xi=a(t)E\neq0$. Therefore the state satisfies Definition~\ref{def:state} but violates macroscopic Prandtl--Reuss normality. For the same data, the stationary fields $u=0$, $\TT=a_0E$, $\theta=1$ and $\Xi=0$, with vanishing production measures, define another relaxed entropy--energy state.
\end{proof}

\begin{remark}[Abstract lift identities do not imply macroscopic normality]
\label{rem:abstract-lift-insufficient}
The preceding example also shows that the limiting lift identities, taken without their generating convergences, do not imply macroscopic normality. Indeed, for the state \eqref{eq:smooth-nonnormal} one may set
\[
\Lambda=\nu=\chi=0\,,\qquad \widetilde\Gamma=0\,,\qquad \Gamma_{\rm pl}=\mu_{\rm pl}\otimes\delta_{(0,0,E)}\,,
\]
together with $\mu_{\rm proj}=0$, $\mu_{\rm over}=\mu_{\rm pl}$, and take all positive-level lifts and barycentres equal to zero. These objects satisfy the graph, marginal, domination, and weighted Minty relations of Theorem~\ref{thm:main}; for the weighted Minty relation this follows from stress admissibility. Nevertheless, their fibre radius is $\rho=0$, whereas the macroscopic radius satisfies $r(\theta)>7/8$. Thus the abstract lift identities alone do not recover macroscopic normality. This construction does not show that the nonnormal state is generated by the coupled approximation. In contrast, with the natural constant regularization of the data, the stationary state above solves every approximate problem.
\end{remark}

\subsection{Complete collapse: rigidity and sharpness}

The fully collapsed regime provides a setting in which the energy and entropy inequalities determine the evolution completely. It also shows that two natural strengthenings of the limiting formulation---measure regularity of the full plastic flow and preservation of the full homogeneous Dirichlet trace---cannot be required even within the approximation-generated class.

\noindent
Let $L^2_\sigma(\Omega)$ denote the $L^2(\Omega;\R^3)$ closure of the smooth compactly supported divergence-free fields.

\begin{proposition}[Rigidity in the fully collapsed case]
\label{prop:collapsed-rigidity}
Assume $\C$ is the identity operator and $r\equiv0$. Let $f=0$, $u_0=0$, $\TT_0=0$, $\theta_0=\bar\theta>0$ be constant, and $u_1=w\in L^2_\sigma(\Omega)$. Then every relaxed entropy--energy state in the sense of Definition~\ref{def:state} satisfies
\begin{equation}\label{eq:collapsed-rigid-state}
\begin{gathered}
u(t)=tw\,,\qquad u_t=w\,,\qquad \TT=0\,,\qquad
\tau=\log\bar\theta\,,\qquad \Theta(t)=\bar\theta\dd x\,,\\
\Xi=\eps(w)\quad\text{in }\mathcal D'(Q;\Sdev)\,,\qquad
\mu_{\rm pl}\restriction_{[0,T)\times\overline\Omega}
=\mu_{\rm th}\restriction_{[0,T)\times\overline\Omega}=0\,.
\end{gathered}
\end{equation}
\end{proposition}

\begin{proof}
Since $r\equiv0$, admissibility gives $\dev\TT=0$. Testing the momentum balance with products of scalar time functions and smooth compactly supported solenoidal fields eliminates both the hydrostatic stress and the thermal term. By density and the initial velocity condition,
\[
\mathbb P_\sigma u_t(t)=w\qquad\text{for a.e. }t\in(0,T)\,,
\]
where $\mathbb P_\sigma$ is the orthogonal projection onto $L^2_\sigma(\Omega)$. Consequently,
\[
\|u_t(t)\|_2^2\ge\|w\|_2^2\,.
\]

On the thermal side, \eqref{eq:state-integrated-entropy} and positivity of the production measures yield
\[
\int_\Omega\tau(t)\dd x\ge|\Omega|\log\bar\theta\,.
\]
The thermal splitting and Jensen's inequality then give
\[
\Theta(t)(\overline\Omega)\ge\int_\Omega e^{\tau(t)}\dd x
\ge|\Omega|\exp\left(\frac1{|\Omega|}\int_\Omega\tau(t)\dd x\right)
\ge|\Omega|\bar\theta\,.
\]
On the other hand, \eqref{eq:state-energy-ineq} yields
\[
\frac12\|u_t(t)\|_2^2+\frac12\|\TT(t)\|_2^2+\Theta(t)(\overline\Omega)
\le\frac12\|w\|_2^2+|\Omega|\bar\theta\,.
\]
Hence all the preceding lower bounds are saturated. Orthogonality gives $u_t=w$, the elastic energy gives $\TT=0$, equality in Jensen's inequality gives $\tau=\log\bar\theta$, and equality of the thermal masses gives $\lambda_t=0$. Since $u(0)=0$, integration in time gives $u(t)=tw$. The constitutive balance then yields
\[
\Xi=\eps(w)\,.
\]
Finally, equality in \eqref{eq:state-integrated-entropy} and positivity imply
\[
\mu_{\rm pl}\restriction_{[0,T)\times\overline\Omega}
=\mu_{\rm th}\restriction_{[0,T)\times\overline\Omega}=0\,.
\]
\end{proof}

\noindent
Proposition~\ref{prop:collapsed-rigidity} applies in particular to every approximation-generated state of Theorem~\ref{thm:main}. The conclusion for the production measures is restricted to $[0,T)$: possible atoms at $t=T$ are not detected by the entropy test class or by the energy inequalities asserted for $t<T$.

\begin{corollary}[The full plastic flow need not be a measure]
\label{cor:nonmeasure-flow}
There are data satisfying Assumption~\ref{ass:standing} for which every approximation-generated state has $\Xi\notin\M_{\rm loc}(Q;\Sdev)$.
\end{corollary}

\begin{proof}
Choose a rectangular box compactly contained in $\Omega$, using translated coordinates with $x_2=0$ in its interior. Let $g\in W^{1,\infty}_c(\R)$ be supported in its second coordinate interval and agree near zero with
\[
g(s)=
\begin{cases}
\displaystyle\int_0^s\sin(1/z)\dd z\,,&s>0\,,\\
0\,,&s\le0\,.
\end{cases}
\]
Then $g'$ is bounded but has infinite variation near zero, so its distributional derivative $g''$ is not a locally finite measure there. Choose smooth compactly supported cutoffs $\eta_1,\eta_3$ in the other two coordinate intervals, equal to one on smaller intervals, and set
\[
\psi(x)=\eta_1(x_1)\eta_3(x_3)g(x_2)\,,\qquad
w=(\partial_2\psi,-\partial_1\psi,0)\,.
\]
The field $w$ is bounded, compactly supported in $\Omega$, and divergence-free. Its mollifications show that $w\in L^2_\sigma(\Omega)$. On a smaller box where the two cutoffs equal one,
\[
\eps_{12}(w)=\frac12g''(x_2)\,.
\]
If $\eps(w)$ were locally a measure, testing with products of a fixed smooth function of $(x_1,x_3)$ of integral one and arbitrary test functions in $x_2$ would imply that $g''$ is a locally finite measure, a contradiction. Proposition~\ref{prop:collapsed-rigidity} therefore gives
\[
\Xi=\eps(w)\notin\M_{\rm loc}(Q;\Sdev)\,.
\]
The distribution $\Xi$ still belongs to $H^{-1}(Q;\Sdev)$, since it is a spatial derivative of an $L^2(Q)$ field. Existence of an approximation-generated state follows from Theorem~\ref{thm:main}.
\end{proof}

\begin{corollary}[Loss of the full Dirichlet trace]
\label{cor:tangential-trace}
There are data satisfying Assumption~\ref{ass:standing} for which approximation-generated states do not retain the full homogeneous Dirichlet condition of the approximate displacements.
\end{corollary}

\begin{proof}
Take $\Omega=B_1(0)$ and $w(x)=b\times x$ for a fixed $b\in\R^3\setminus\{0\}$. This smooth field is divergence-free, has zero normal trace, and has a nonzero tangential trace. Smooth radial cutoffs $\chi_j$ with compact support in the ball and tending to one give
\[
\Div(\chi_jw)=\nabla\chi_j\cdot w=0\,,\qquad \chi_jw\longrightarrow w\quad\text{in }L^2(\Omega;\R^3)\,,
\]
so $w\in L^2_\sigma(\Omega)$. Proposition~\ref{prop:collapsed-rigidity} forces
\[
u(t)=tw\,,
\]
whose full boundary trace is nonzero for $t>0$. Here $\eps(w)=0$, so this example is independent of the loss of measure regularity established in Corollary~\ref{cor:nonmeasure-flow}.
\end{proof}

\noindent
The preceding results use the admissible fully collapsed law $r\equiv0$. They show that, under the standing assumptions, neither measure regularity of the full plastic flow nor preservation of the full homogeneous Dirichlet trace can be required in general. They do not imply the same loss for every nontrivial thermally truncated law or on regions where the yield radius remains uniformly positive.

\subsection{Limits of macroscopic identification}

There is also a purely geometric obstruction to recovering macroscopic normality by averaging the lifted variables. Let $C_0=\overline B_1^{\Sdev}(0)$ and let $e\in\Sdev$ be a unit deviator. Both fibre states $(Z,\omega)=(e,e)$ and $(Z,\omega)=(-e,-e)$ satisfy
\[
Z\in C_0\,,\qquad \ip{Z}{\omega}=h_0(\omega)\,.
\]
Their mixture with weights $3/4$ and $1/4$ has stress and direction barycentres $\overline Z=\frac12e$ and $\overline\omega=\frac12e$. Since $|\overline Z|<1$, the stress barycentre lies strictly inside $C_0$, where $N_{C_0}(\overline Z)=\{0\}$, whereas $\overline\omega\neq0$. Thus support of a lift on the exact normality graph does not by itself imply a normality relation for the corresponding barycentres, even when $C_0$ is strictly convex.

\noindent
For comparison, suppose that on a given region a plastic measure $P$ and a compatible stress--flow pairing are available and satisfy, in the von Mises case,
\begin{equation}\label{eq:conditional-macro-work}
[\dev\TT:P]=r(e^\tau)|P|\,,\qquad [\dev\TT:P]\le|\dev\TT|\,|P|\,.
\end{equation}
If, for some $\eta>0$,
\[
|\dev\TT|\le r(e^\tau)-\eta \qquad |P|\text{-a.e.}\,,
\]
then
\[
r(e^\tau)|P| =[\dev\TT:P] \le|\dev\TT|\,|P| \le\bigl(r(e^\tau)-\eta\bigr)|P|\,.
\]
Hence $\eta|P|\le0$ and therefore $P=0$ on that region. This is the usual exclusion of plastic flow in the strict elastic interior. Proposition~\ref{prop:global-lifts} provides its exact geometric counterpart on the lifted fibres, where
\[
\ip{Z}{\omega} =\rho h_0(\omega) =\max_{S\in\rho C_0}\ip{S}{\omega}\,.
\]
It does not, however, yield the macroscopic relation \eqref{eq:conditional-macro-work} for the distribution $\Xi$.

\noindent
The preceding observations do not weaken the constitutive conclusions of Section~\ref{sec:plastic}. Rather, they identify their precise scope. Exact normality and maximum plastic work are retained on the lifted fibres, while a corresponding relation for the macroscopic pair $(\TT,\Xi)$ requires additional identification or duality that is not available under the standing assumptions.

\section{Concluding discussion}
\label{sec:conclusion}

Complete thermal collapse preserves a closed instantaneous normality graph but removes uniform coercivity of the unweighted plastic rate. Logarithmic entropy and total energy control the thermal coupling and the distance to the moving constraint, while support weighting and compact-fibre lifting retain finite constitutive measures together with exact fibrewise normality and maximum plastic work. The positive-radius lift makes the nondegenerate restrictions compatible at all continuity levels. On fixed open regions satisfying \eqref{eq:fixed-region-positive}, the selected restrictions at levels $0<\eta\le\delta$ identify the full distributional plastic flow as a finite measure. The resulting theorem constructs approximation-generated relaxed entropy--energy states for the stated finite-energy data, while the available compactness does not identify a canonical macroscopic stress--plastic-flow pairing.

The results of Section~\ref{sec:normality} clarify the sharpness of this formulation. The intrinsic axioms alone define a strictly broader balance class, as shown by the smooth nonnormal state, whereas the fully collapsed rigidity result shows that neither global measure regularity of the full plastic flow nor preservation of the full homogeneous Dirichlet trace can be required under the standing assumptions, even within the approximation-generated class. These limitations do not affect the lifted constitutive relations established in Section~\ref{sec:plastic}; rather, they distinguish the fibrewise normality retained by the approximation from additional macroscopic identifications that would require further structure or a compatible stress--plastic-flow duality. The geometric appendix records the constraint estimate and graph compactness that remain independent of the homothetic form of the yield sets.

\appendix
\section{A geometric principle for moving convex sets}
\label{app:geometry}

The constitutive compactness argument of Section~\ref{sec:plastic} uses only a small part of the homothetic structure of the yield sets. We record here the underlying convex-geometric mechanism in a form that applies to general Hausdorff-continuous families of compact convex sets. This result concerns only the constraint estimate and the lifted normality structure; the thermomechanical existence argument remains tied to the assumptions imposed in the main body of the paper.

\noindent
Let $E$ be a finite-dimensional Euclidean space, let $\mathcal A$ be a compact metric space, and let $\alpha\mapsto C(\alpha)$ be a Hausdorff-continuous family of nonempty compact convex subsets of $E$ satisfying
\[
0\in C(\alpha)\subset\overline B_R(0) \qquad\text{for every }\alpha\in\mathcal A\,.
\]
We denote the corresponding support function by $h(\alpha,\omega)=\sup_{Z\in C(\alpha)}\ip{Z}{\omega}$. Let $Y\subset\R^m$ be bounded and open, and let $\alpha_n:Y\to\mathcal A$ and $S_n:Y\to E$ be measurable. Define
\[
Z_n=\Pi_{C(\alpha_n)}S_n\,,\qquad d_n=|S_n-Z_n|\,,\qquad \omega_n=\frac{S_n-Z_n}{d_n}\,,\qquad \xi^n=a_n(d_n)\omega_n
\]
on $\{d_n>0\}$. On $\{d_n=0\}$ we set $\xi^n=0$ and choose $\omega_n$ arbitrarily with $|\omega_n|=1$. Finally, set $q^n=\ip{\xi^n}{S_n}$.
\begin{proposition}[Shape-independent overstress control and lifted normality]
\label{prop:abstract}
Assume $\sup_n\int_Y\sqrt{q^n}\dd y<\infty$. Then the estimates \eqref{eq:algebraic-rate}, \eqref{eq:main-rate}, and \eqref{eq:overstress-tail} remain valid with $Q$ replaced by $Y$. Let $\delta>0$ and let $B_{n,\delta}\subset Y$ be measurable sets such that
\[
\overline B_\delta(0)\subset C(\alpha_n(y)) \qquad\text{for every }y\in B_{n,\delta}\,.
\]
If
\[
\sup_n\int_{B_{n,\delta}} h(\alpha_n,\omega_n)|\xi^n|\dd y<\infty\,,
\]
then $\int_{B_{n,\delta}}|\xi^n|\dd y\le\frac{C}{\delta}$. Moreover, every weak-star limit of
\[
(\operatorname{id},\alpha_n,Z_n,\omega_n)_\# \bigl(\mathbf1_{B_{n,\delta}}|\xi^n|\dd y\bigr)
\]
in $\M^+\!\left( \overline Y\times\mathcal A\times \overline B_R(0)\times\{\omega\in E:|\omega|=1\} \right)$ is concentrated on the closed set
\begin{equation}\label{eq:abstract-graph}
Z\in C(\alpha)\,,\qquad \ip{Z}{\omega}=h(\alpha,\omega)\,.
\end{equation}
\end{proposition}
\begin{proof}
By the characterization of the metric projection onto a closed convex set, $\omega_n\in N_{C(\alpha_n)}(Z_n)$ on $\{d_n>0\}$. Hence $\ip{Z_n}{\omega_n}=h(\alpha_n,\omega_n)$ and therefore
\[
q^n = a_n(d_n)\bigl(h(\alpha_n,\omega_n)+d_n\bigr) \ge a_n(d_n)d_n\,,
\]
because $0\in C(\alpha_n)$. Taking square roots gives exactly the scalar bound used in Proposition~\ref{prop:overstress}. The proof of that proposition is therefore unchanged and yields \eqref{eq:algebraic-rate}, \eqref{eq:main-rate}, and \eqref{eq:overstress-tail} with $Q$ replaced by $Y$. On $B_{n,\delta}$, the inclusion $\overline B_\delta(0)\subset C(\alpha_n(y))$ gives $h(\alpha_n,\omega_n)\ge\delta$ because $|\omega_n|=1$. Consequently,
\[
\delta\int_{B_{n,\delta}}|\xi^n|\dd y \le \int_{B_{n,\delta}} h(\alpha_n,\omega_n)|\xi^n|\dd y \le C\,,
\]
which proves the asserted mass bound and hence weak-star compactness of the associated lifts. It remains to identify their support. Hausdorff continuity of $\alpha\mapsto C(\alpha)$ implies that $\{(\alpha,Z)\in\mathcal A\times\overline B_R(0):Z\in C(\alpha)\}$ is closed. It also implies joint continuity of the support function $(\alpha,\omega)\longmapsto h(\alpha,\omega)$ on the compact set $\mathcal A\times\{\omega\in E:|\omega|=1\}$. Therefore the set defined by \eqref{eq:abstract-graph} is closed. Since every approximate lift is concentrated on this set, the same is true for each of its weak-star limits.
\end{proof}

\noindent
Proposition~\ref{prop:abstract} isolates the convex-geometric part of the constitutive compactness argument and shows that neither the overstress estimate nor the closedness of the lifted normality graph requires a homothetic family of yield sets. It does not extend the thermomechanical existence theorem itself: such an extension would additionally require a compatible fixed-level stress construction and an entropy estimate controlling the corresponding support work.

\section*{Statements and Declarations}

\noindent\textbf{Funding.}
The author received no specific funding for this work.

\smallskip
\noindent\textbf{Competing interests.}
The author has no relevant financial or non-financial interests to disclose.

\smallskip
\noindent\textbf{Data availability.}
No datasets were generated or analysed during the current study.

\smallskip
\noindent\textbf{Use of AI-assisted tools.}
ChatGPT (OpenAI) was used during manuscript preparation for language and LaTeX editing and for literature checking.

\end{document}